\documentclass[a4paper,12pt]{article}
\pdfoutput=1
\usepackage{epsfig,amsmath,amsfonts,amsthm}
\usepackage[figuresright]{rotating}
\usepackage[left=28mm,right=20mm,top=30mm,bottom=20mm]{geometry}
\usepackage{subfig}
\usepackage{graphicx}
\usepackage{rotating}
\usepackage{booktabs}
\usepackage{moreverb}

\DeclareGraphicsExtensions{.pdf,.png,.jpg}
\usepackage{natbib}
\numberwithin{equation}{section}

\def\e{\hbox{E}}

\def\min{\hbox{min}}
\def\max{\hbox{max}}

\newtheorem{theorem}{Theorem}[section]
\newtheorem{lemma}[theorem]{Lemma}

\providecommand{\keywords}[1]{\textbf{{Keywords}} #1}
\usepackage[nokeyprefix]{refstyle}
\usepackage{varioref}
\usepackage{xr}
\usepackage{hyperref}
\newcommand{\thickbar}[1]{\mathbf{\bar{\text{$#1$}}}}

\begin{document}

\title{Point and interval estimators of a changepoint in  stochastical dominance between two distributions}

\author{Elena Kulinskaya and  David C. Hoaglin }

\date{\today}

\maketitle

\begin{abstract}
For differences between means of continuous data from independent groups, the customary scale-free measure of effect is the standardized mean difference (SMD).  To justify use of SMD, one should be reasonably confident that the group-level variances are equal. Empirical evidence often contradicts this assumption. Thus, we have investigated an alternate approach, based on stochastic ordering of the treatment and control distributions, that takes into account means and variances. For applying stochastic ordering, our development yields a key quantity, $\mathsf{A}$, the outcome value at which the direction of the ordering of the treatment and control distributions changes. 

Using an extensive simulation, we studied relative bias of point estimators of $\mathsf{A}$ and coverage and relative width of bootstrap confidence intervals.

\end{abstract}


\section{Introduction}

The simplicity of a scale-free effect size for comparisons such as the difference between a treatment mean and a control mean, $\mu_T - \mu_C$, has led to widespread use of a standard deviation, $\sigma$, as the denominator. The resulting measure of effect size is $\delta = (\mu_T - \mu_C) / \sigma$, a standardized mean difference (SMD).

The within-group variances, $\sigma_T^2$ and $\sigma_C^2$, provide the basis for $\sigma$. When interest focuses on the effect of the treatment relative to an untreated control group,  \cite{Glass1981} (page 29) take $\sigma = \sigma_C$.
 When the studies do not involve an untreated control group (e.g., the studies used small and large classes in estimating effects of class size), \cite{Cohen_1965} and \cite{Glass1981} (page 38) assume $\sigma^2 = \sigma_T^2 = \sigma_C^2$. (Various authors have discussed and modified this form of SMD.) Actual meta-analyses work with the sample means, $\bar{X}_T$ and $\bar{X}_C$, and variances, $s_T^2$ and $s_C^2$. The assumption that $\sigma_T^2 = \sigma_C^2$ leads to the pooled estimate of $\sigma^2$, $s_p^2 = [(n_T - 1)s_T^2 + (n_C - 1)s_C^2] / (n_T + n_C - 2)$. It also permits translation of $\delta$ into measures of nonoverlap between the treatment and control distributions \citep{Cohen1988}.

In practice, $\sigma_T^2$ and $\sigma_C^2$ are seldom equal. \cite{Grissom2001} list various applications in which one expects $\sigma_T^2 > \sigma_C^2$ and others in which one expects $\sigma_T^2 < \sigma_C^2$. When $\sigma_T^2 \neq \sigma_C^2$, two complications arise. First, $(\bar{X}_T - \bar{X}_C) / s_p$, even after correction for bias \citep{hedges1983random}, no longer estimates a well-defined $\delta$. Second, no agreement exists on a definition of $\sigma$ based on $\sigma_T^2$  and $\sigma_C^2$.

For a fuller picture, we draw on the notion of stochastic ordering. Our approach gives a description of the relation between treatment and control that is more informative than $\mu_T - \mu_C$ and avoids the limitation $\sigma_T = \sigma_C$.

Section~\ref{sec:stochastic} discusses stochastic ordering and derives $\mathsf{A}$, the outcome value at which the direction of the ordering of the treatment and control distributions changes. Section~\ref{sec:estA} develops estimates of $\mathsf{A}$. Section~\ref{sec:sims} discusses the design and results of our simulations for point and interval estimation of $\mathsf{A}$. The figures in Appendices A to C show the full results of our simulations.

\section{Stochastic ordering} \label{sec:stochastic}

For two random variables $X_1$ and $X_2$ from a family of continuous distributions on $(-\infty, \infty)$, $X_1$ is stochastically less than $X_2$ if $X_2$ tends to have larger values than $X_1$. In terms of probabilities, the definition is
\begin{equation} \label{eq:FSD}
P(X_1 > A) \leq P(X_2 > A) \text{   for all   } -\infty < A < \infty.
\end{equation}
This relation, denoted by $X_1\preceq X_2$, is known as first-order stochastic dominance of $X_1$ by $X_2$ (\cite{marsha112011}). (\cite{Lehmann_1975} discusses the role of stochastic ordering in tests based on ranks.)
Such complete stochastic ordering is often too restrictive. A more realistic version is stochastic ordering on an interval $\mathcal{A}$ (i.e., for all $A \in \mathcal{A}$), which we denote by $X_1 \preceq _{\mathcal{A}} X_2$ (i.e., $X_2$ stochastically dominates $X_1$ on $\mathcal{A}$). That is, $ P(X_2 > A) - P(X_1 > A) \geq 0 \text{   for all   } A \in {\mathcal{A}}$. On the complement of $\mathcal{A}$, denoted by $\thickbar{\mathcal{A}}$, $X_2$ is stochastically less than $X_1$: $X_2 \preceq _{\thickbar{\mathcal{A}}} X_1$. Ordinarily, $\mathcal{A}$ is half-open, either $\mathcal{A} = [A_{min}, \infty)$ or $\mathcal{A} = (-\infty, A_{max}]$.

To apply stochastic dominance (on an interval) in comparing treatment with control, we relate the endpoints $A_{min}$ and $A_{max}$ to $\mu_T$, $\mu_C$, $\sigma_T$, and $\sigma_C$. We assume that $X_1$ and $X_2$ have cumulative distribution functions $F_i(x) = F(x; \mu_i, \sigma_i)$ belonging to a continuous location-scale family on $(-\infty, \infty)$: $F_i(x) = F(\frac{x - \mu_i} {\sigma_i})$. The following Lemma provides the main results.
\begin{lemma} \label{lemma}
Let $X_i \sim F(x;\mu_i, \sigma_i)$, $i=1,2$.
\begin{enumerate}
\item When $\sigma_1 = \sigma_2$, $X_1 \preceq X_2$ if and only if $\mu_1 \leq \mu_2$.
\item When $\sigma_1 > \sigma_2$, $X_1 \preceq_{\mathcal{A}} X_2$  and $X_2 \preceq_{\thickbar {\mathcal{A}}} X_1$ for ${\mathcal{A}} = \{A \leq \mathsf{A} = \mu_2 + (\mu_1 - \mu_2) / (1 - \sigma_1 / \sigma_2) \}$.
\item When $\sigma_1 < \sigma_2$, $X_1 \preceq_{\mathcal{A}} X_2$ and $X_2 \preceq_{\thickbar {\mathcal{A}}} X_1$ for ${\mathcal{A}} = \{A \geq \mathsf{A} = \mu_2 + (\mu_1 - \mu_2) / (1 - \sigma_1 / \sigma_2) \}$.
\end{enumerate}
\end{lemma}
{\bf Proof}
Because $P(X_i > A) = 1 - F_i(A)$, Equation~(\ref{eq:FSD}) is equivalent to $F(\frac{A - \mu_2} {\sigma_2}) \leq F(\frac{A - \mu_1} {\sigma_1})$. Because $F$ is an increasing function,
$\frac{A - \mu_2} {\sigma_2} \leq \frac{A - \mu_1} {\sigma_1}$, which can be rewritten as
$$A \leq \mu_2 + (\mu_1 - \mu_2) / (1 - \sigma_1 / \sigma_2).$$
The three statements follow directly.

Statement 1 covers the case $\sigma_1 = \sigma_2$. In Statements 2 and 3 we see that the endpoints that we denoted by $A_{max}$ and $A_{min}$ are equal; we denote that shared endpoint by $\mathsf{A}$. Also, $ P(X_2 > A) = P(X_1 > A) $ when $A = \mathsf{A}$.

For a simple example, we let $\mu_C = 1$, $\sigma_C = 1$, $\mu_T = 1$, and $\sigma_T = 2$. Then $\mathsf{A} = \mu_T + (\mu_C - \mu_T) / (1 - \sigma_C / \sigma_T) = \mu_T = 1$. Because $\sigma_C < \sigma_T$, we apply Statement 3: $X_C \preceq_{\mathcal{A}} X_T$ and $X_T \preceq_{\thickbar {\mathcal{A}}} X_C$ for ${\mathcal{A}} = \{A \geq 1 \}$. What we infer about the comparison between treatment and control depends on the direction of the outcome.  If larger values of $X$ are more favorable, treatment is better than control: larger values (above 1) are more likely under the treatment than under the control. For the same reason, if larger values of $X$ are less favorable, control is better than treatment. These results are consistent with the values of $\mu_T$, $\mu_C$, $\sigma_T$, and $\sigma_C$: because $\mu_T = \mu_C$ and $\sigma_T > \sigma_C$, the probability for $X_T$ is spread out more, both below and above $\mu_T = \mu_C = 1$. Thus, stochastic ordering yields more information than SMD = 0, which suggests no difference.

In a comparative study, let control correspond to $i = 1$ and treatment to $i = 2$. If lower values of the continuous outcome are favorable (and $\sigma_1 > \sigma_2$ or $\sigma_1 = \sigma_2$), then a successful treatment would be stochastically less than control over a suitably large interval of values of $A$ (i.e., for a large subset of patients).

In applications, we estimate  $\mathsf{A}$, using $\bar{X}_T$, $\bar{X}_C$, $s_T$, and $s_C$ (in place of $\mu_T$, $\mu_C$, $\sigma_T$, and $\sigma_C)$, and must accept the resulting uncertainty.

\section{Estimation of $\mathsf{A}$} \label{sec:estA}

Let $X_i \sim N(\mu_i, \sigma_i)$ for $i = 1,2$. In Statements 2 and 3 of Lemma 1, the endpoint of the interval $\mathcal{A}$ is $\mathsf{A} = \mu_2 + (\mu_1 - \mu_2) / (1 - \sigma_1 / \sigma_2)$. Without loss of generality, let $\sigma_1 < \sigma_2$. (If $\sigma_1 > \sigma_2$, the derivation below holds with subscripts 1 and 2 interchanged.)

The factor $\Sigma = 1 / (1 - \sigma_1/\sigma_2)$ is the sum of the geometric series $1 + (\sigma_1 / \sigma_2) + (\sigma_1 / \sigma_2)^2 + \cdots $. The plug-in estimator
\begin{equation} \label{estA}
\hat{\mathsf{A}} = \bar{X}_2 + (\bar{X}_1 - \bar{X}_2) / (1 - s_1/s_2).
\end{equation}
has $\e(\hat{\mathsf{A}}) = \mu_2 + (\mu_1 - \mu_2) \e[(1 - s_1/s_2)^{-1}]$ (the sample mean and variance of normal random variables are independent), but $\e[(1 - s_1/s_2)^{-1}] \neq (1 - \sigma_1 / \sigma_2)^{-1}$. In fact, this expectation is not finite, because of variability from $s_2$. (Similarly, the $k$th moment of the $F$ distribution $F_{\nu_1, \nu_2}$ is finite only when $k < \nu_{2} / 2$.)  However, the individual $(\sigma_1 / \sigma_2)^j$ in the geometric series do have finite expectation for $j < n_2 - 1$. Thus, we can estimate the corresponding $(s_1 / s_2)^j$.

 As $(s_1 / s_2)^2 \sim (\sigma_1^2 / \sigma_2^2) F_{n_1 - 1, n_2 - 1}$, the moments of $s_1 / s_2$ around zero are given by $K_j = \e[(s_1 / s_2)^j] = (\sigma_1 / \sigma_2)^j M_{j/2} (F_{n_1 - 1, n_2 - 1})$, for $j < n_2 - 1$, where
 \begin{equation} \label{Mj}
 M_{j/2} (F_{n_1 - 1, n_2 - 1}) = \left( \frac{n_2 - 1} {n_1 - 1} \right)^{j/2} \frac {\Gamma((n_1 + j - 1)/2) \Gamma((n_2 - j - 1)/2)} {\Gamma((n_1 - 1)/2) \Gamma((n_2 - 1)/2)}.
 \end{equation}
 Thus, we estimate $(\sigma_1 / \sigma_2)^j$ for $j < n_2 - 1$ by
 \[ (s_1 / s_2)^j / M_{j/2} (F_{n_1 - 1, n_2 - 1}). \]
Then, for $k < n_2 - 1$, we define
\begin{equation} \label{approx2}
\hat\Sigma_k = \sum _0^k M_{j/2}^{-1} (F_{n_1 - 1, n_2 - 1}) (s_1 / s_2)^j
\end{equation}
and
\begin{equation} \label{approx}
\hat{\mathsf{A}}_k = \bar{X}_2 + (\bar{X}_1 - \bar{X}_2) \hat{\Sigma}_k.
\end{equation}
If the contributions of $(s_1 / s_2)^j$ decrease rapidly enough as $j$ increases, a value of $k$ substantially less than $n_2 - 2$ may suffice. (We examine this possibility in Section~\ref{sec:sims}.)

The second moment  of $\hat{\mathsf{A}}_k$ includes  an expectation of
$$\sum _{j=0}^k\sum _{h=0}^k M_{j/2}^{-1} (F_{n_1 - 1, n_2 - 1})M_{h/2}^{-1} (F_{n_1 - 1, n_2 - 1}) (s_1/s_2)^{j+h}.$$ Thus, the variance of $\hat{\mathsf{A}}_k$ exists for $2k<n_2-1$.  However, the distribution of $\hat{\mathsf{A}}_k$ is skewed, and the  variance is of no practical importance.

To find a confidence interval for $\mathsf{A}$, we use a parametric bootstrap.

\section{Simulations for point and interval estimation of $\mathsf{A}$} \label{sec:sims}

\subsection{Design} \label{sec:DesignSims}

Our simulations pursued two objectives. First, we wanted to evaluate the speed of convergence of the estimators $\hat{\mathsf{A}}_k$ in $k$. For that,
only estimation of $\Sigma = (1 - \sigma_1 / \sigma_2)^{-1}$ by $\hat\Sigma_k$ in Equation~(\ref{approx2}) is of interest. Therefore, in one series of simulations, we generated 10,000 sample values of $\hat{\Sigma} \sim 1/(1 - (\alpha^2 F(n_T - 1, n_C - 1))^{1/2})$ for each combination of a value of $\alpha = \sigma_T / \sigma_C = 0.1(0.1)0.9$ and a value of $n = n_T = n_C = 12,\;22,\;32,\;42,\;52,\;77,\;102,\;152,\;202$.  For each value of $n$, we obtained the values of $\hat\Sigma_k$ for $k\leq k_{max} = n - 2$.

In practice, many studies' group sizes fall in the ranges covered by our choices. Many studies allocate subjects equally to the two groups, and rough equality holds more widely (as in the studies analyzed by \cite{Rubio-Aparicio_2018_BehavResMeth_2057}).
The values of $\alpha$  were chosen  to represent the range of values at both sides of the 5\% critical values of the $F$-test for  equality of variances vs $\alpha < 1$.  For the above group sizes, these  values are 0.596,  0.693, 0.741, 0.771, 0.793, 0.827, 0.848, 0.874 and  0.890, respectively. The corresponding 1\% critical values are 0.473, 0.592, 0.652, 0.691, 0.719, 0.764, 0.792, 0.827 and  0.848.

In a second series of simulations we evaluated the bias of $\hat{\mathsf{A}}_k$ and the coverage and length of parametric-bootstrap-based confidence intervals for ${\mathsf{A}}$. In the first series of simulations, calculation of $\hat\Sigma_k$  required relatively little time, so the second series used only $k = k_{max} = \min(n - 2, 500)$.

We varied four parameters: the difference between the means ($\Delta = \mu_T - \mu_C$), the standard deviation in the control arm ($\sigma_C$), the ratio of the standard deviations ($\alpha = \sigma_T / \sigma_C$), and the sample size in each group ($n$). We kept the mean in the control group constant at $\mu_C = 1$ as it does not seem to matter (cf. Equation~(\ref{approx})).  Table~\ref{tab:design1} lists the values of each parameter. We generated $M = 1000$ repetitions for each combination of parameters, and each repetition used $B = 2000$ bootstrap samples to obtain confidence intervals at the 95\% and 90\% confidence levels. We also obtained the bootstrap median value.

In repetition $j$, we calculated $\mu_{T} = \mu_C + \Delta$ and $\sigma_{T} = \alpha \sigma_{C}$, generated  sample means from a normal distribution $\bar{X}_{ij} \sim N(\mu_{i}, \sigma_{i}^2 / n_{i})$ and sample variances from a scaled chi-square distribution $s_{ij}^2 \sim \sigma_{i}^2 \chi^2_{n_{i} - 1} / (n_{i} - 1)$. We used Equation~(\ref{approx2}) and Equation~(\ref{approx}) to calculate $\hat{\mathsf{A}}_{kj}$ for $k = k_{max}$. From these we calculated relative bias of $\hat{\mathsf{A}}_{kj}$ as $(\hat{\mathsf{A}}_{kj} - {\mathsf{A}}) / {\mathsf{A}} \times 100\%$.

Next, in each repetition, we calculated $B$ bootstrap samples of the sample means and variances, using the generated sample means $\bar X_{ij}$ and sample variances $s_{ij}^2$ as the true means and variances, and obtained the estimated values $\hat{\mathsf{A}}_{kj(b)}$. From these we obtained sample quantiles at $p = 0.025, \;0.05,\;0.5, \;0.95$ and $0.975$ to evaluate coverage of 95\% and 90\% confidence intervals for ${\mathsf{A}}$ and their relative width on the left and right of ${\mathsf{A}}$.

Overall, we considered 600 combinations of parameters. The simulations used \emph{R} statistical software \citep{rrr}. User-friendly R programs implementing our methods are available at  https://osf.io/ag7vm .

\begin{table}[ht]
	\caption{ \label{tab:design1} \emph{Values of parameters in the second series of simulations}}
	\begin{footnotesize}
		\begin{centering}
			\begin{tabular}
				{|l|l|}
				\hline
				Parameter & Values \\
                                	& \\
				\hline
				$n$ (size of each arm) & 20, 50, 100, 250, 500, 1000 \\
				\hline
                			$\mu_C$ (mean in the control arm) & 1\\
                			\hline		
                			$\Delta$ (difference between the means, $\Delta = \mu_T - \mu_C$) & $-2$, $-1$,  1, 2 \\
               			\hline
                			$\sigma_C$ (standard deviation in the control arm)& 0.2, 0.5, 1, 1.5, 2 \\
                			\hline                			
				$\alpha$ (ratio of standard deviations, $\alpha = \sigma_T / \sigma_C$) & $0.3$, $0.5$, $0.7$, $0.9$ \\
				\hline
			\end{tabular}
		\end{centering}
	\end{footnotesize}
\end{table}

\subsection{Results} \label{sec:ResultSims1}

\begin{figure}[ht]
	\centering
	\includegraphics[scale=0.35]{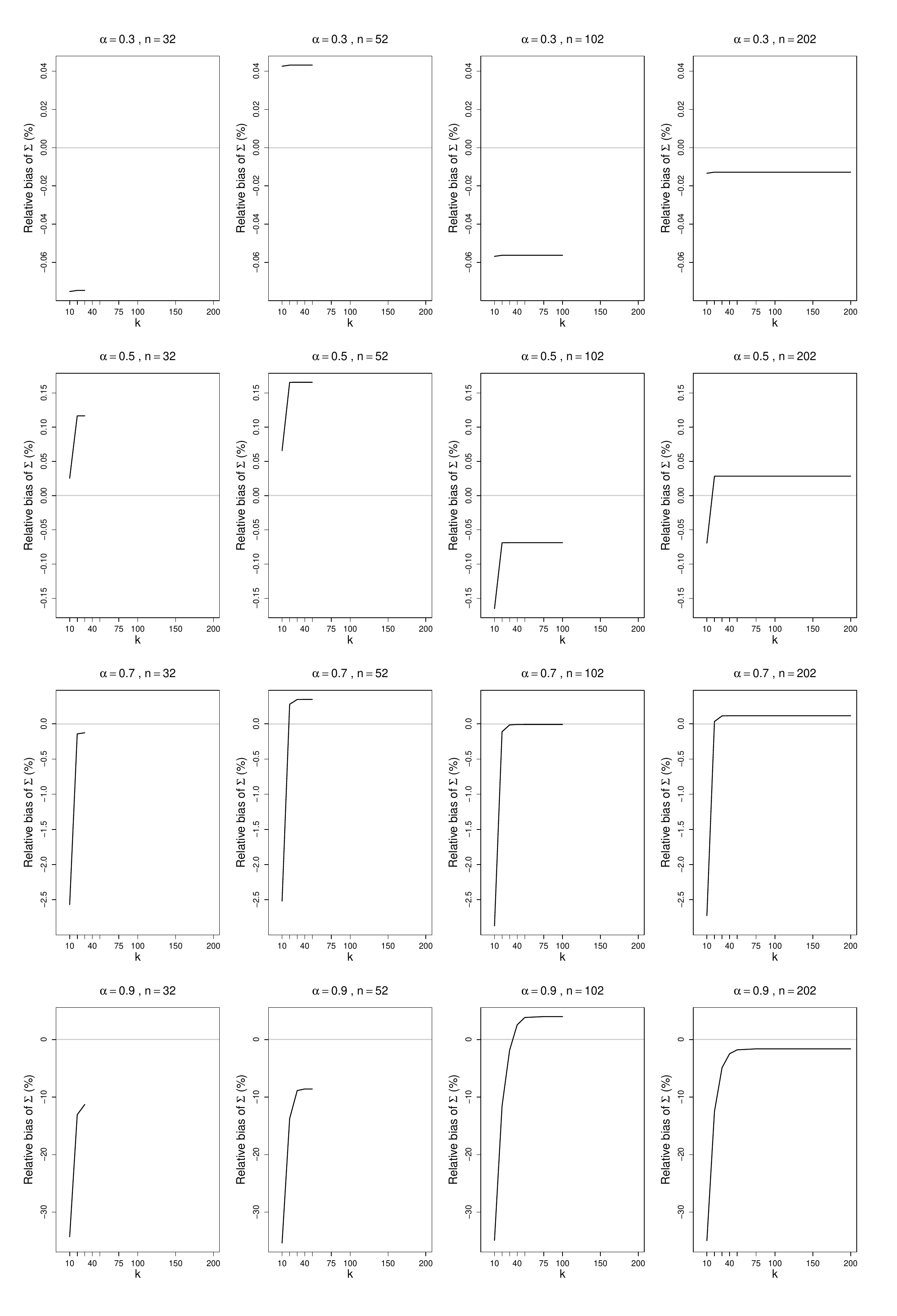}
	\caption{Relative bias of $\hat{\Sigma}_k$ for $k \leq k_{max} = n - 2$ as a function of $k = 10,\;20,\;30,\;40,\;50,\;75,\;100,\;150,\;200$
when $n = n_C = n_T = 32,\;52,\;102,\;202$ and $\alpha = 0.3,\;0.5,\;0.7,\;0.9$. Light grey line at zero.}
   	\label{Fig0}
\end{figure}

\subsubsection{Convergence of the estimators $\hat{\Sigma}_k$ in $k$ (Figure~\ref{Fig0})}
Bias in estimating $\Sigma$ by the plug-in estimator  $\hat\Sigma=(1-s1/s2)^{-1}$  is very considerable for large values of $\alpha$. For $\alpha = 0.7$, relative bias is 30\% for $n = 102$, and for $\alpha = 0.9$ relative bias is 281\% for $n = 102$ and $-220$\% for $n = 202$.

%%%> bias
%%%          [,1]        [,2]        [,3]        [,4]
%%%[1,]   1.353836   0.8772781    0.348746   0.1881394
%%%[2,]   7.321042   3.6998555    1.530691   0.8042377
%%%[3,]  11.371541  30.3885888    8.793579   3.9513328
%%%[4,] -79.542708 281.0312463 -219.529493 -32.3017842

Bias of the estimators  $\hat\Sigma_k$ of $\Sigma$ may be considerable for small values of $k$, especially when $\alpha=0.9$. Estimation  is fast and quickly converges in $k$. Bias appears to be low and stable for $k \geq 50$ (Figure~\ref{Fig0}).  As time is not a limiting factor, we recommend taking $k = k_{max}$ for estimation of $\mathsf{A}$.

\subsubsection{Relative bias of the estimators $\hat{\mathsf{A}}_k$ (Appendix A)}

Bias of $\hat{\mathsf{A}}_k$ for $k = k_{max}$ does not seem to depend much on $\Delta$ and $\sigma_C$. However, it increases dramatically when  the ratio of standard deviations $\alpha$ approaches 1, and slowly decreases in the  sample size $n$.

For $n = 50$, relative bias of  $\hat{\mathsf{A}}_k$ is within 2\% for $\alpha \leq 0.5$, within 7\% for $\alpha = 0.7$ and above 77\% for $\alpha = 0.9$.  When $n = 100$, relative bias is within 2.5\% for $\alpha = 0.7$ but still reaches 57\% for $\alpha = 0.9$. However, this bias decreases to 26\% for $n = 250$, to 8\% for $n = 500$ and to 2\% for $n = 1000$.

Bootstrap median values are typically lower than those of $\hat{\mathsf{A}}_k$ as the distribution is skewed to the right.

\subsubsection{Coverage of the bootstrap confidence intervals of  $\mathsf{A}$ (Appendix B)}

Coverage of ${\mathsf{A}}$ is reasonable for moderate $n$ when the ratio of standard deviations $\alpha \leq 0.7$, but deteriorates dramatically when $\alpha$ approaches 1, and slowly improves as the sample size $n$ increases.

For $n = 50$, coverage is within 2 percentage points of nominal 95\% for $\alpha \leq 0.7$, but drops to 10\% for $\alpha = 0.9$ when $\sigma_C = 0.3$ and then increases in $\sigma_C$, reaching  about 50\% when $\sigma_C = 2$ and $|\Delta| = 1$. It reaches about 24\% when $|\Delta| = 2$.
However, coverage does improve as $n$ increases and does not seem to depend on $\Delta$ and $\sigma_C$ when $n \geq 100$. For $\alpha = 0.9$, coverage at the 95\% nominal level is about 80\% at 95\% nominal level, increases to 92\% for $n = 250$, to 93--94\% for $n = 500$ and to 95--96\% for $n = 1000$.

\subsubsection{Relative width of the bootstrap confidence intervals for $\mathsf{A}$ (Appendix C)}

%The relative length of the left and right half-widths of the bootstrap $(1 - \beta)$-confidence intervals measures relative distance from ${\mathsf{A}}$ to the bootstrap quantiles at the $\beta / 2$ and $(1 - \beta/2)$ levels. This provides information on the total width of the confidence intervals, and also on the skewness of the respective distributions.

%The width of the intervals increases with $\sigma_C$ and  $\alpha$, and slowly decreases in $n$. It also depends on $\Delta$.

%For $n = 50$, the left half-width is greater than to the right half-width (i.e., the distributions are skewed to the left), for $\alpha \leq 0.5$, but skewed to the right for $\alpha \geq 0.7$ (Figure~\ref{Fig5}).  For $100 \leq n \leq 500$, the distribution is skewed to the left for $\alpha\leq 0.7$, but skewed to the right for $\alpha = 0.9$ (not shown). When $n = 1000$, the distribution is skewed to the left even for $\alpha = 0.9$ for negative $\Delta$, but is still skewed to the right for $\Delta \geq 1$. The intervals are considerably narrower for large $n$.

To examine the width of the $100 (1 - \beta)\%$ confidence interval for $\mathsf{A}$, and also the skewness of the bootstrap distributions of $\hat{\mathsf{A}}$, we express the distance from the $100 \beta / 2$ bootstrap percentile ($\hat{\mathsf{A}}_{\beta / 2}^{\ast}$) to $\mathsf{A}$ (known in the simulations) and the distance from $\mathsf{A}$ to the $100 (1 - \beta / 2)$ percentile ($\hat{\mathsf{A}}_{1 - \beta / 2}^{\ast}$) as percentages of $|\mathsf{A}|$. Specifically,
\[ \hat{W}_{left} = 100 (\mathsf{A} - \hat{\mathsf{A}}_{\beta / 2}^{\ast}) / | \mathsf{A} | \]
and
\[ \hat{W}_{right} = 100 (\hat{\mathsf{A}}_{1 - \beta / 2}^{\ast} - \mathsf{A}) / | \mathsf{A} | . \]
We are interested in the size of these relative widths, and $\mathsf{A}$ can be negative. Using percentages of $ | \mathsf{A} | $ makes them more nearly comparable across variation in $\mathsf{A}$, which is a function of $\mu_C$, $\Delta$, and $\alpha$. For example, when $\mu_C = 1$ and $\alpha$ and $\Delta$ have the values used in the simulation, $\mathsf{A}$ ranges from $-1.86$ to $+3.86$ when $\alpha = 0.3$ and from $-19.0$ to $+21.0$ when $\alpha = 0.9$.

The sum of $\hat{W}_{left}$ and $\hat{W}_{right}$ equals the width of the confidence interval, relative to $  | \mathsf{A} | $, and the difference between $\hat{W}_{left}$ and $\hat{W}_{right}$ reflects the skewness of the bootstrap distribution of $\hat{\mathsf{A}}$.

The figures in Appendix C plot, versus $\sigma_C$, $\hat{W}_{left}$ and $\hat{W}_{right}$ for the 95\% and 90\% confidence intervals for $\mu_C = 1$ and various combinations of $n$, $\alpha$, and $\Delta$. In each plot, $\mathsf{A}$ varies with $\sigma_C$ only via $\alpha$, but $\sigma_C$ and $\sigma_T = \alpha \sigma_C$ both contribute to variation in the samples produced in the parametric bootstrap.

In Figure~\ref{FigC2} (for $n = 50$), the relative width tends to increase with $\sigma_C$, most noticeably when $\Delta = -1$, and it is greater when $\Delta < 0$ than when $\Delta > 0$. Skewness has a more complicated pattern: when $\alpha \leq 0.5$, $\hat{W}_{left} > \hat{W}_{right}$ for $\Delta < 0$, but $\hat{W}_{right} > \hat{W}_{left}$ for $\Delta > 0$; when $\alpha \geq 0.7$, that pattern holds with the signs of $\Delta$ reversed. Also, the difference between $\hat{W}_{left}$ and $\hat{W}_{right}$ tends to increase with $\alpha$.

For other $n$, the same summary of total (relative) width applies. For $100 \leq n \leq 500$, the sign change in the pattern of skewness occurs between $\alpha = 0.7$ and $\alpha = 0.9$, and for $n = 1000$ the pattern remains the same. As it did for $n = 50$, the difference between $\hat{W}_{left}$ and $\hat{W}_{right}$ tends to increase with $\alpha$.

%\section*{Acknowledgments}
%The work by E. Kulinskaya was supported by the Economic and Social Research Council
%[grant number ES/L011859/1].

\clearpage
\bibliography{Refs2DH.bib}
\clearpage

%\end{document}

\section*{Appendices}
\begin{itemize}
\item Appendix A: Relative bias of the estimators $\hat{\mathsf{A}}_k$
\item Appendix B: Coverage of the bootstrap confidence intervals for $\mathsf{A}$
\item Appendix C: Relative width of the bootstrap confidence intervals for $\mathsf{A}$

\end{itemize}

\setcounter{figure}{0}
\setcounter{section}{0}
\clearpage

\section*{Appendix A: Plots of relative bias of the estimators $\hat{\mathsf{A}}_k$}

Each figure corresponds to $\mu_C = 1$ and a value of the within-group sample size $n$ (= 20, 50, 100, 250, 500, 1000).\\

\noindent For each combination of a value of $\alpha = \sigma_T / \sigma_C$ (= 0.3, 0.5, 0.7, 0.9) and a value of the difference in means $\Delta = \mu_T - \mu_C$ (= $-2, -1, 1, 2$), a panel plots, versus the standard deviation in the Control group $\sigma_{C}$  (= $0.2, 0.5, 1, 1.5, 2$), the relative bias of $\hat{\mathsf{A}}_k$ for $k = k_{max} = \max(n - 2, 500)$ (black) and the bootstrap median (green) from the parametric bootstrap with $B = 2000$.

\clearpage

\setcounter{figure}{0}
\setcounter{section}{0}
\renewcommand{\thefigure}{A.\arabic{figure}}

%%%%%%%%%%%%%%%%%%%%%%%%%%%%%%%%%%%%%%%%%%%%%%%%%%%%%%%%%%%%%%%%%%%%%%%%%%%%%%%%%%%%%%%%%%%%%%%%%%%%%%%%%%%%%%%%
%%%%%%%%%%%%%%%%%%%%%%%%%%%%%%%%%%%%%%%% delta_C=-2.5  %%%%%%%%%%%%%%%%%%%%%%%%%%%%%%%%%%%%%%%%%%%%%%%%%%%%%%%%%%

\begin{figure}[ht]
	\centering
	\includegraphics[scale=0.35]{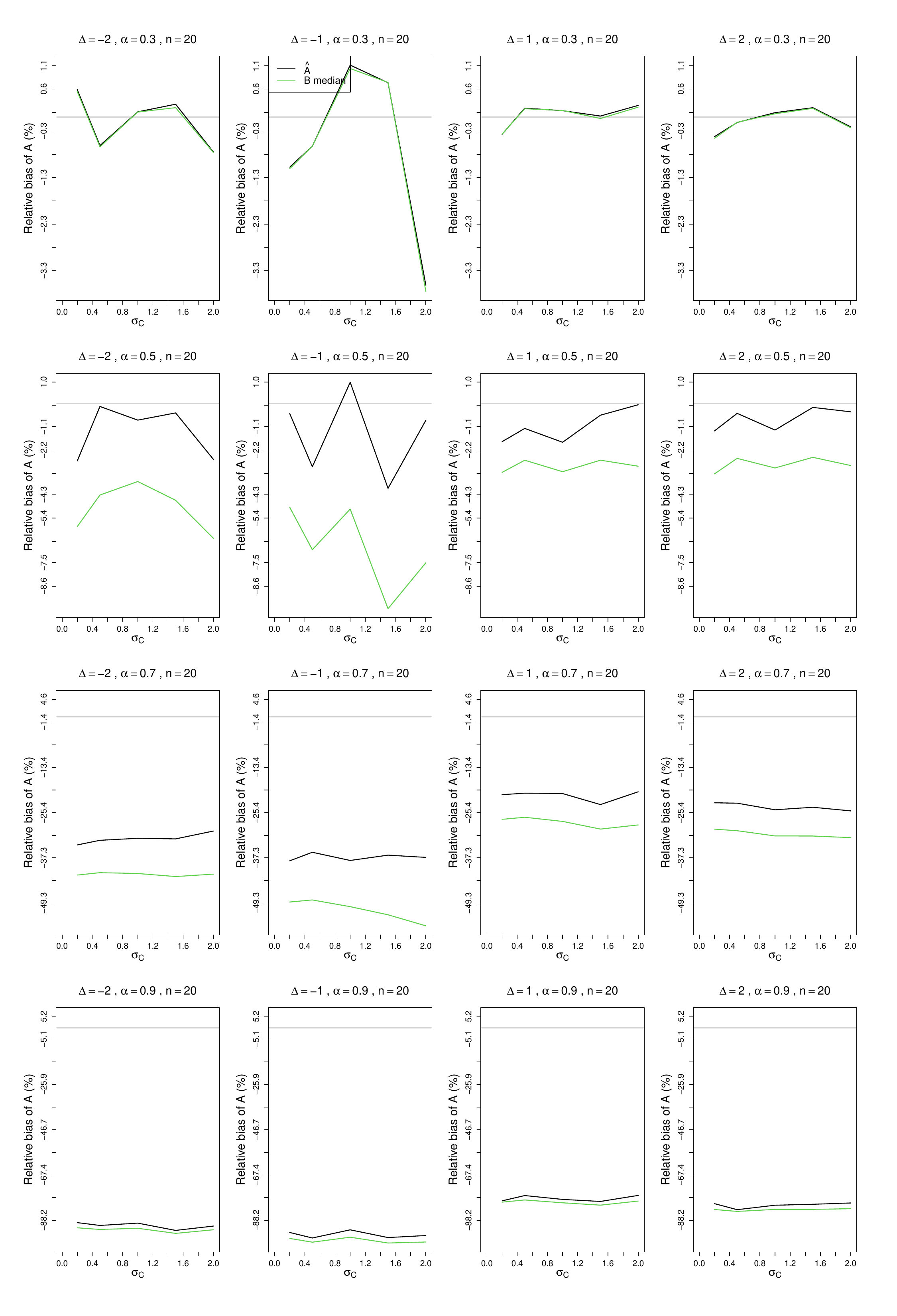}
	\caption{Relative bias of $\hat{\mathsf{A}}_k$ for $k = k_{max}$ (black) and the bootstrap median (green) as a function of the standard deviation in the Control group $\sigma_C$ when $\mu_C = 1$, $n = 20$, $\alpha = 0.3,\;0.5,\;0.7,\;0.9$, and $\Delta = -2,\;-1,\; 1,\;2$. Light grey line at zero.}
   	\label{FigA1}
\end{figure}

\begin{figure}[ht]
	\centering
	\includegraphics[scale=0.35]{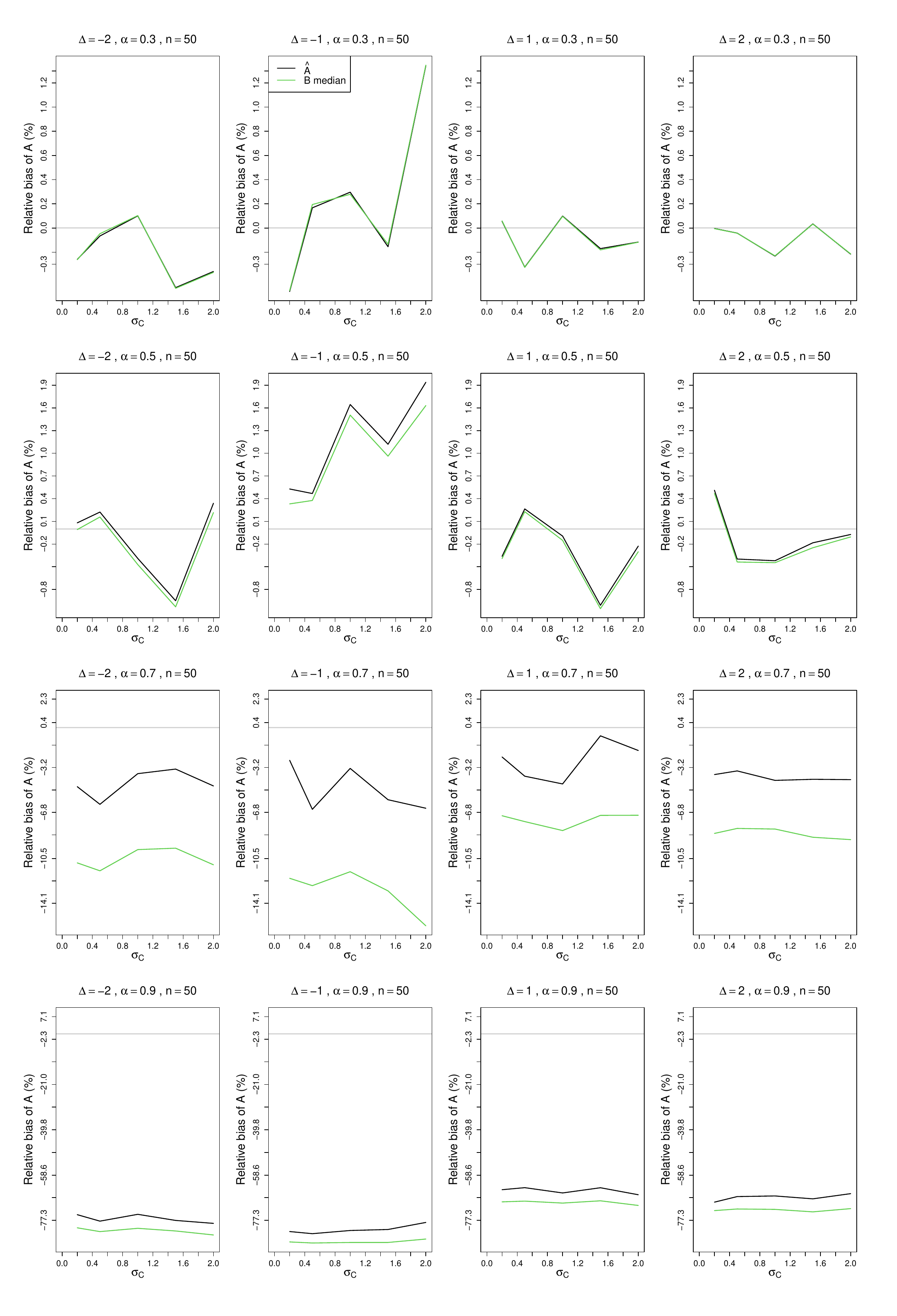}
	\caption{Relative bias of $\hat{\mathsf{A}}_k$ for $k = k_{max}$ (black) and the bootstrap median (green) as a function of the standard deviation in the Control group $\sigma_C$ when $\mu_C = 1$, $n = 50$, $\alpha = 0.3,\;0.5,\;0.7,\;0.9$, and $\Delta = -2,\;-1,\; 1,\;2$. Light grey line at zero.}
   	\label{FigA2}
\end{figure}

\begin{figure}[ht]
	\centering
	\includegraphics[scale=0.35]{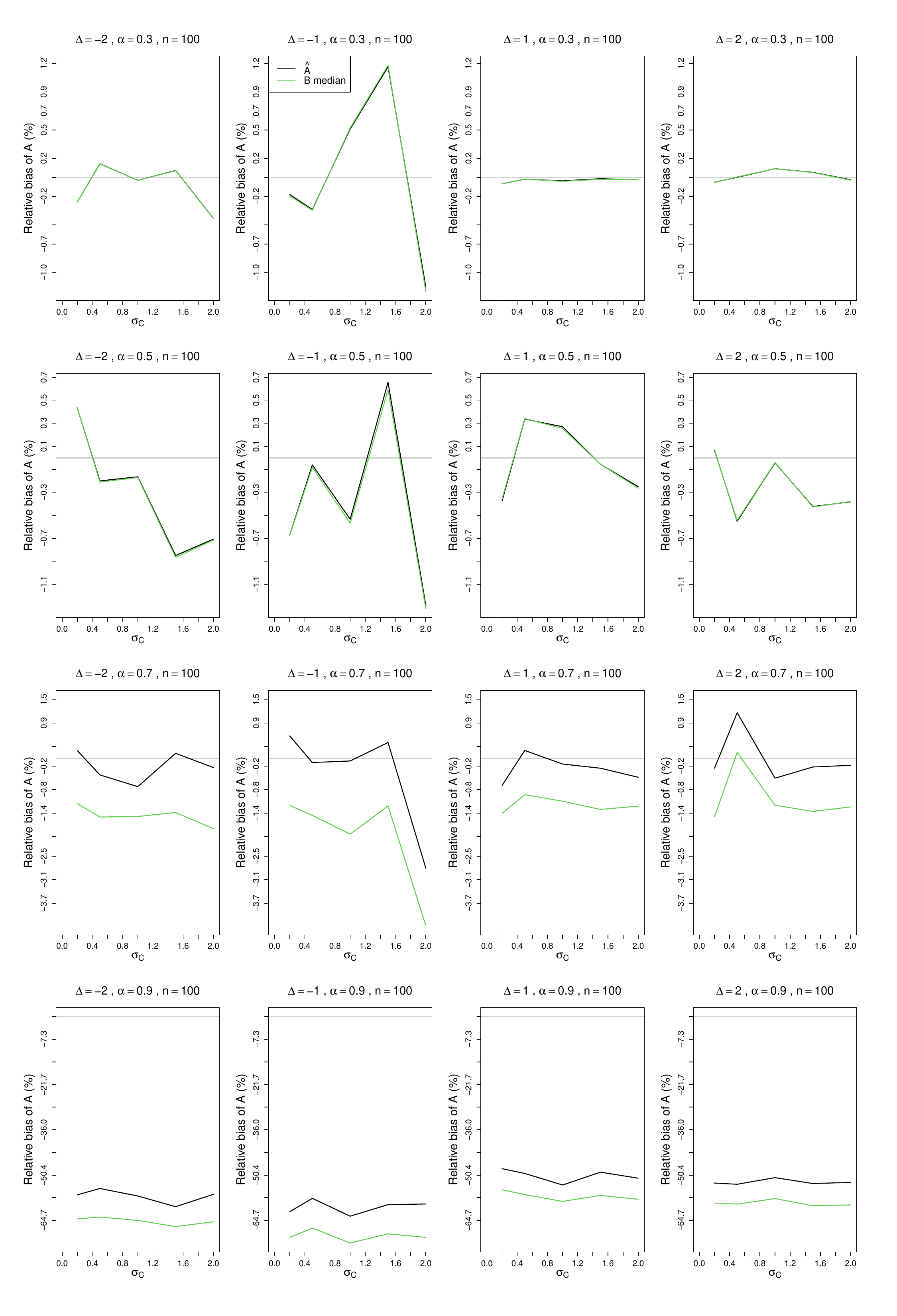}
	\caption{Relative bias of $\hat{\mathsf{A}}_k$ for $k = k_{max}$ (black) and the bootstrap median (green) as a function of the standard deviation in the Control group $\sigma_C$ when $\mu_C = 1$, $n = 100$, $\alpha = 0.3,\;0.5,\;0.7,\;0.9$, and $\Delta = -2,\;-1,\; 1,\;2$. Light grey line at zero.}
   	\label{FigA3}
\end{figure}

\begin{figure}[ht]
	\centering
	\includegraphics[scale=0.35]{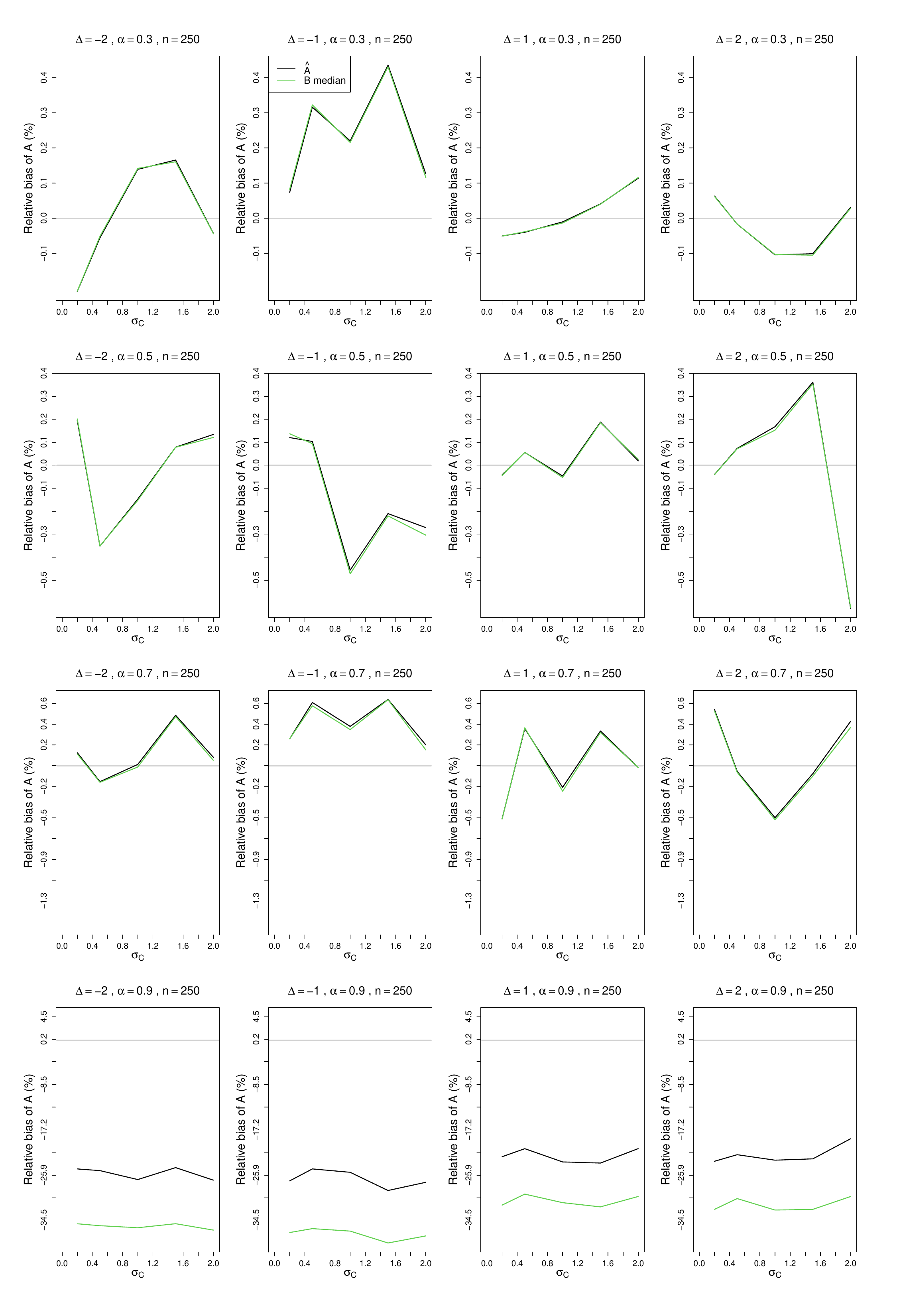}
	\caption{Relative bias of $\hat{\mathsf{A}}_k$ for $k = k_{max}$ (black) and the bootstrap median (green) as a function of the standard deviation in the Control group $\sigma_C$ when $\mu_C = 1$, $n = 250$, $\alpha = 0.3,\;0.5,\;0.7,\;0.9$, and $\Delta = -2,\;-1,\; 1,\;2$. Light grey line at zero.}
   	\label{FigA4}
\end{figure}

\begin{figure}[ht]
	\centering
	\includegraphics[scale=0.35]{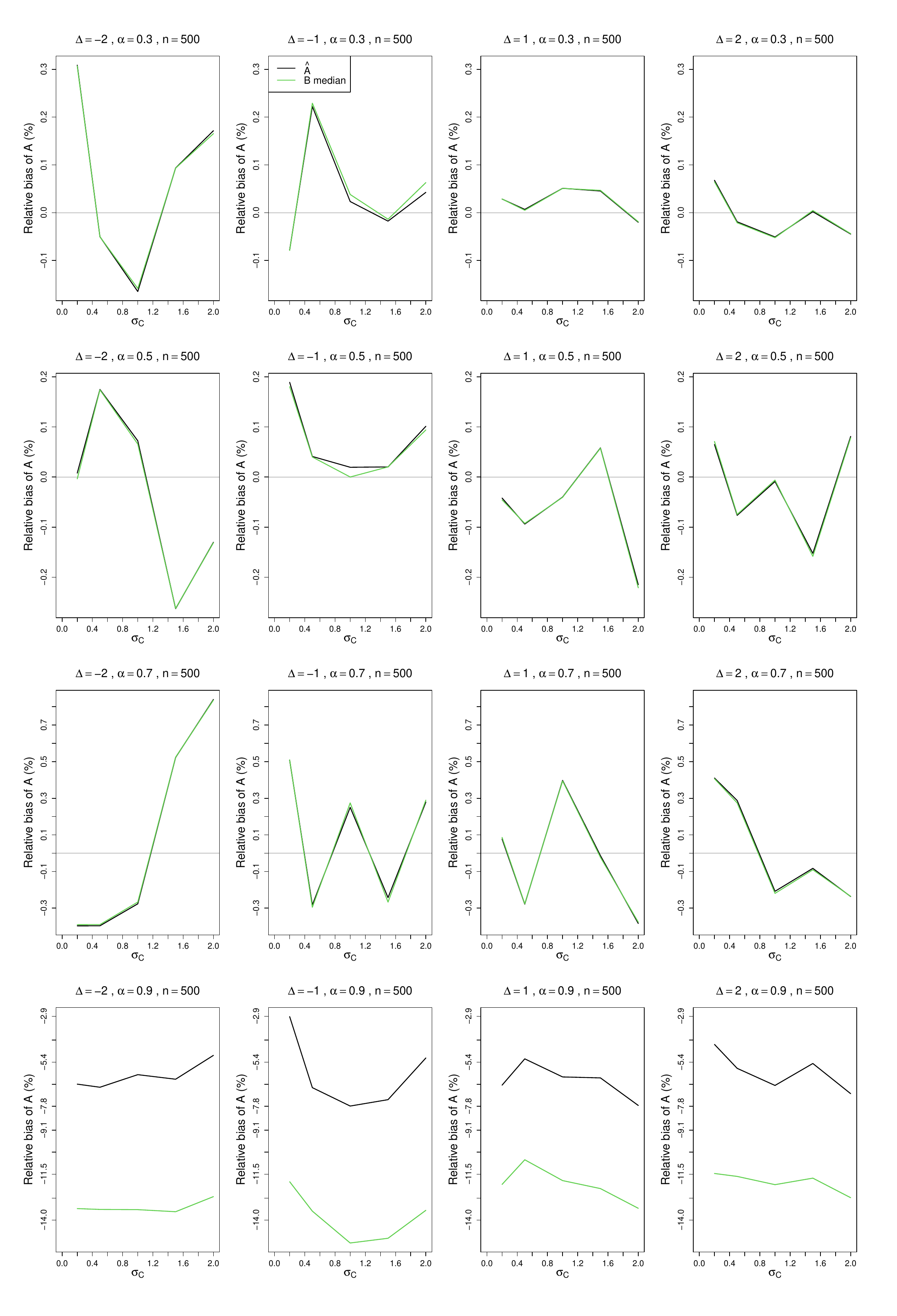}
	\caption{Relative bias of $\hat{\mathsf{A}}_k$ for $k = k_{max}$ (black) and the bootstrap median (green) as a function of the standard deviation in the Control group $\sigma_C$ when $\mu_C = 1$, $n = 500$, $\alpha = 0.3,\;0.5,\;0.7,\;0.9$, and $\Delta = -2,\;-1,\; 1,\;2$. Light grey line at zero.}
   	\label{FigA5}
\end{figure}

\begin{figure}[ht]
	\centering
	\includegraphics[scale=0.35]{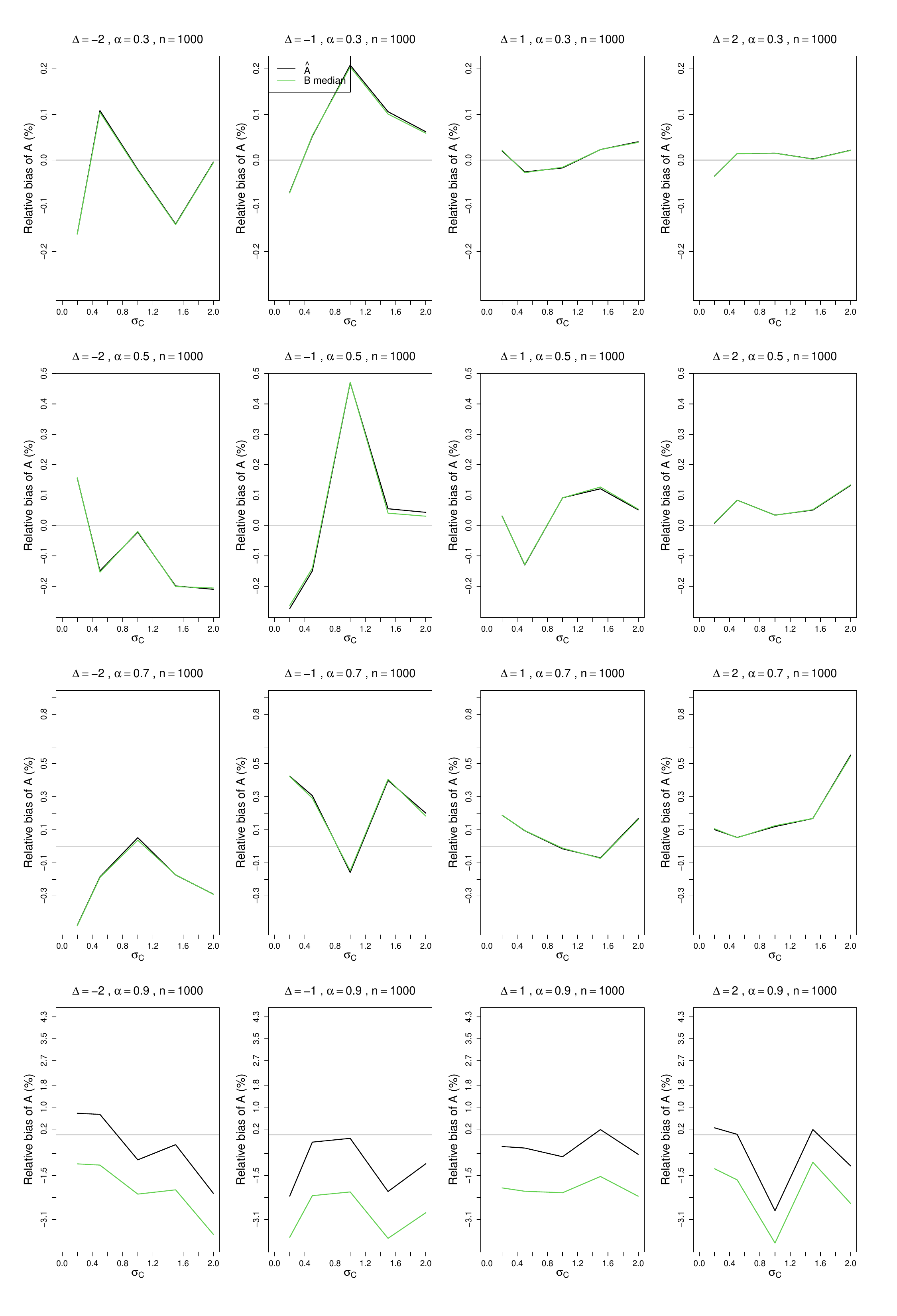}
	\caption{Relative bias of $\hat{\mathsf{A}}_k$ for $k = k_{max}$ (black) and the bootstrap median (green) as a function of the standard deviation in the Control group $\sigma_C$ when $\mu_C = 1$, $n = 1000$, $\alpha = 0.3,\;0.5,\;0.7,\;0.9$, and $\Delta = -2,\;-1,\; 1,\;2$. Light grey line at zero.}
   	\label{FigA6}
\end{figure}

\clearpage

\setcounter{figure}{0}
\setcounter{section}{0}

\section*{Appendix B: Coverage of the bootstrap confidence intervals for $\mathsf{A}$}

Each figure corresponds to $\mu_C = 1$ and a value of the within-group sample size $n$ (= 20, 50, 100, 250, 500, 1000).\\

\noindent For each combination of a value of $\alpha = \sigma_T / \sigma_C$ (= 0.3, 0.5, 0.7, 0.9) and a value of the difference in means $\Delta = \mu_T - \mu_C$ (=$-2, -1, 1, 2$), a panel plots, versus the  standard deviation in the Control group $\sigma_{C}$  (= $0.2, 0.5, 1, 1.5, 2$), coverage of $\mathsf{A}_k$ for $k = k_{max} = \max(n - 2, 500)$ by parametric-bootstrap-based confidence intervals with $B = 2000$ at the 95\% (black) and 90\% (green) levels.

\clearpage
\renewcommand{\thefigure}{B.\arabic{figure}}

\begin{figure}[ht]
	\centering
	\includegraphics[scale=0.35]{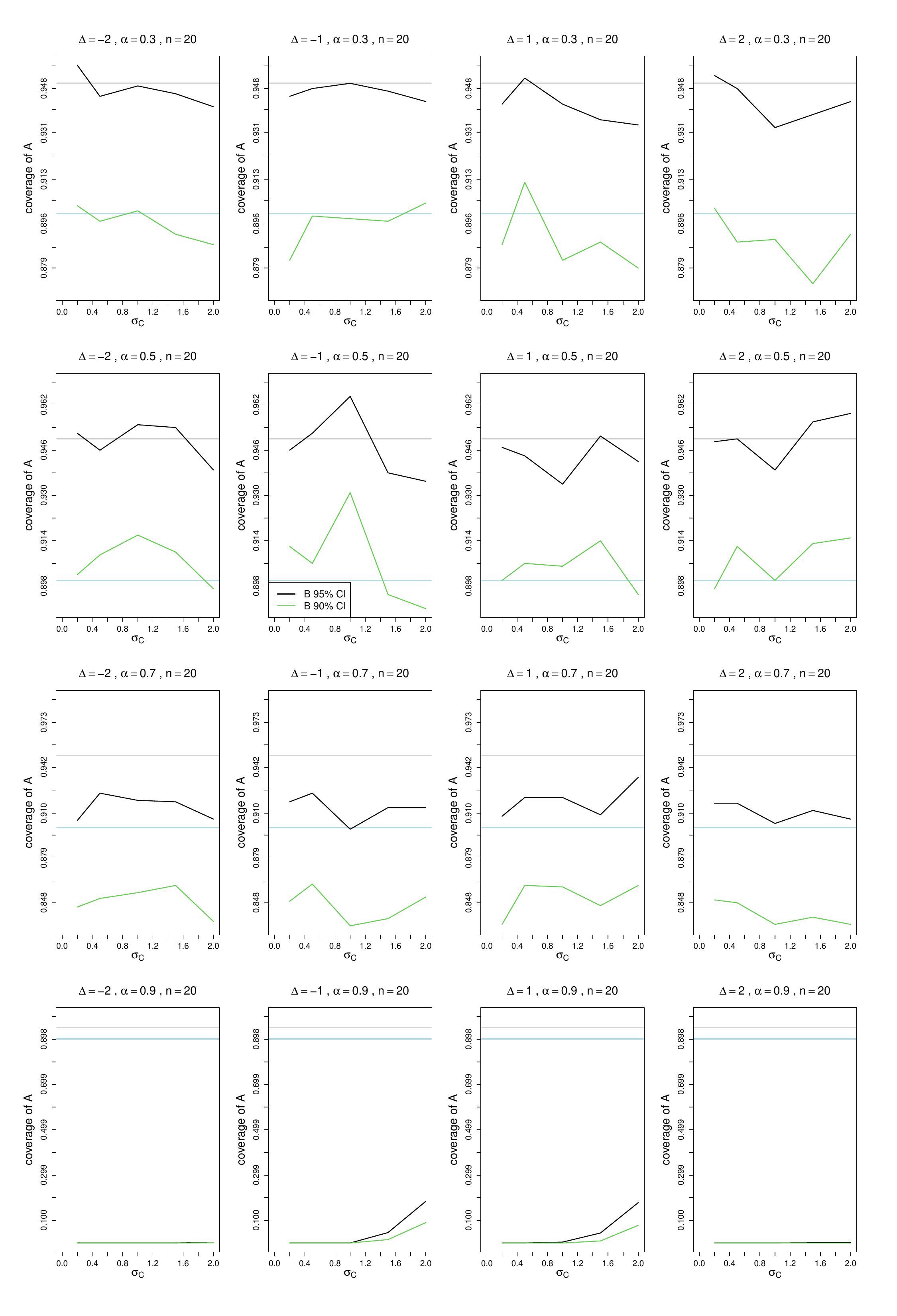}
	\caption{Coverage of $\mathsf{A}$ for $k = k_{max}$ at 95\% (black) and 90\% (green) nominal levels as a function of the standard deviation in the Control group $\sigma_C$ when $\mu_C = 1$, $n = 20$, $\alpha = 0.3,\;0.5,\;0.7,\;0.9$, and $\Delta = -2,\;-1,\; 1,\;2$. Light grey line at .95 and light blue line at .90.}
   	\label{FigB1}
\end{figure}

\begin{figure}[ht]
	\centering
	\includegraphics[scale=0.35]{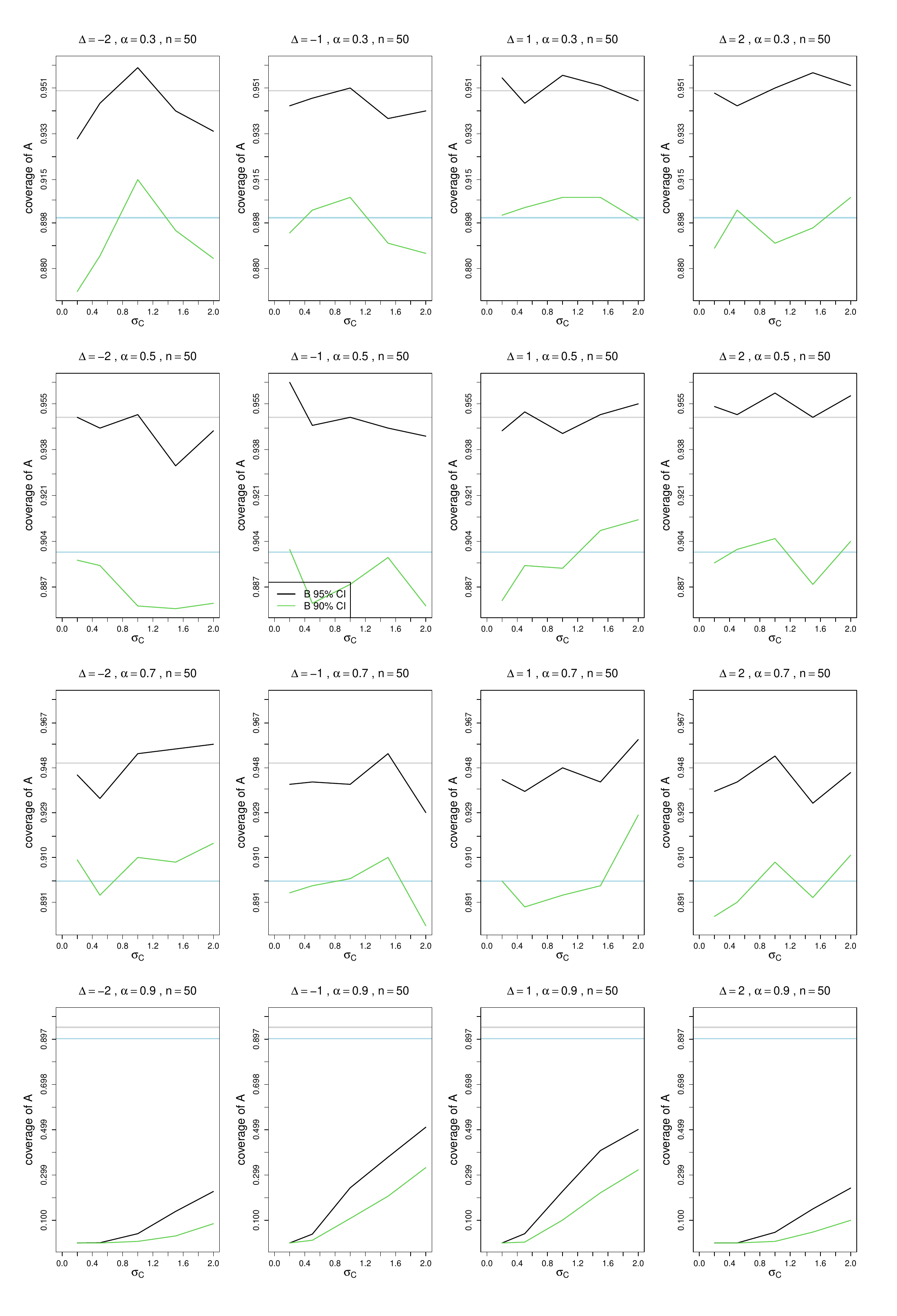}
	\caption{Coverage of $\mathsf{A}$ for $k = k_{max}$ at 95\% (black) and 90\% (green) nominal level as a function of the standard deviation in the Control group $\sigma_C$ when $\mu_C = 1$, $n = 50$, $\alpha = 0.3,\;0.5,\;0.7,\;0.9$, and $\Delta = -2,\;-1,\; 1,\;2$. Light grey line at .95 and light blue line at .90.}
   	\label{FigB2}
\end{figure}

\begin{figure}[ht]
	\centering
	\includegraphics[scale=0.35]{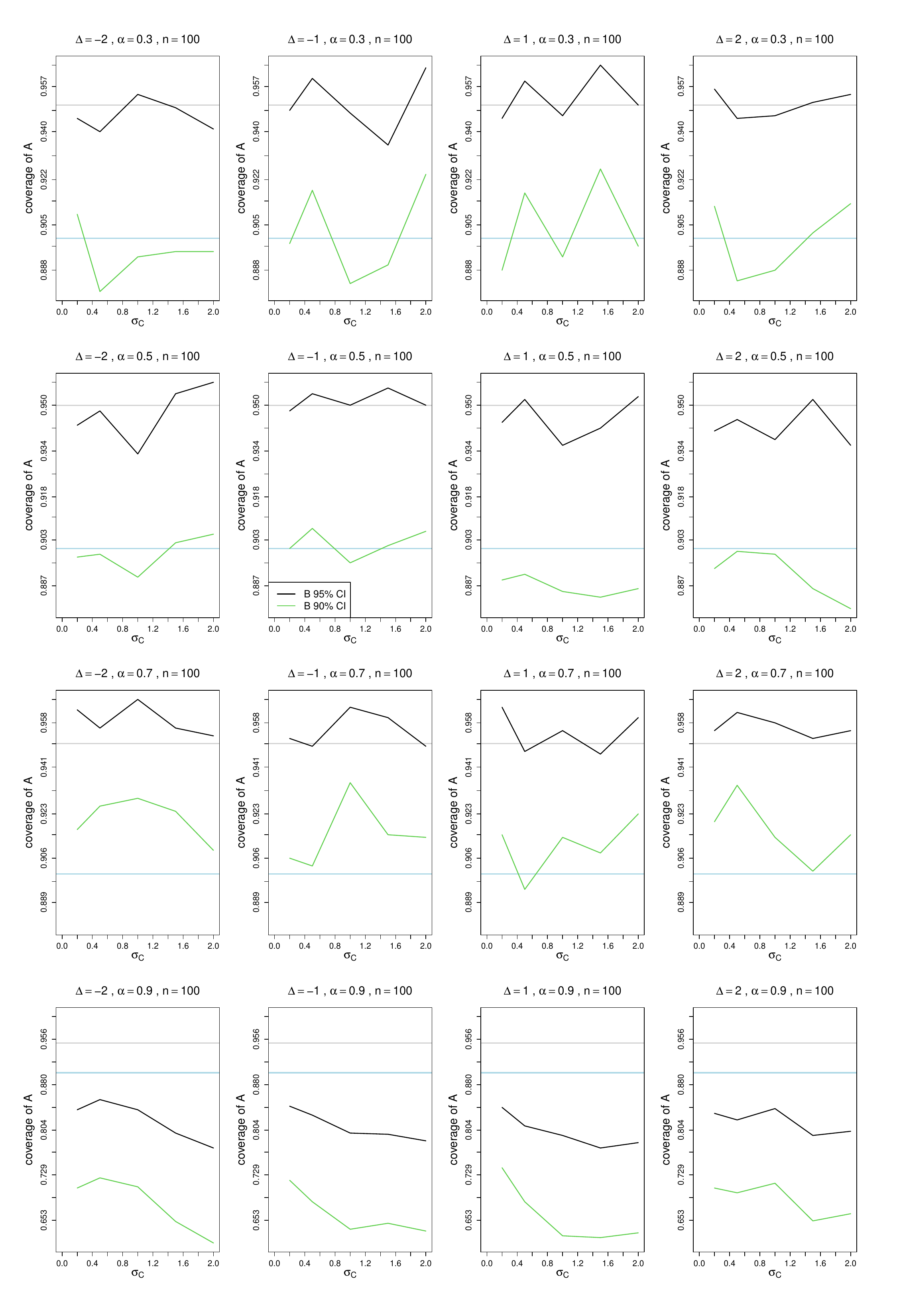}
	\caption{Coverage of $\mathsf{A}$ for $k = k_{max}$ at 95\% (black) and 90\% (green) nominal level as a function of the standard deviation in the Control group $\sigma_C$ when $\mu_C = 1$, $n = 100$, $\alpha = 0.3,\;0.5,\;0.7,\;0.9$, and $\Delta = -2,\;-1,\; 1,\;2$. Light grey line at .95 and light blue line at .90.}
   	\label{FigB3}
\end{figure}

\begin{figure}[ht]
	\centering
	\includegraphics[scale=0.35]{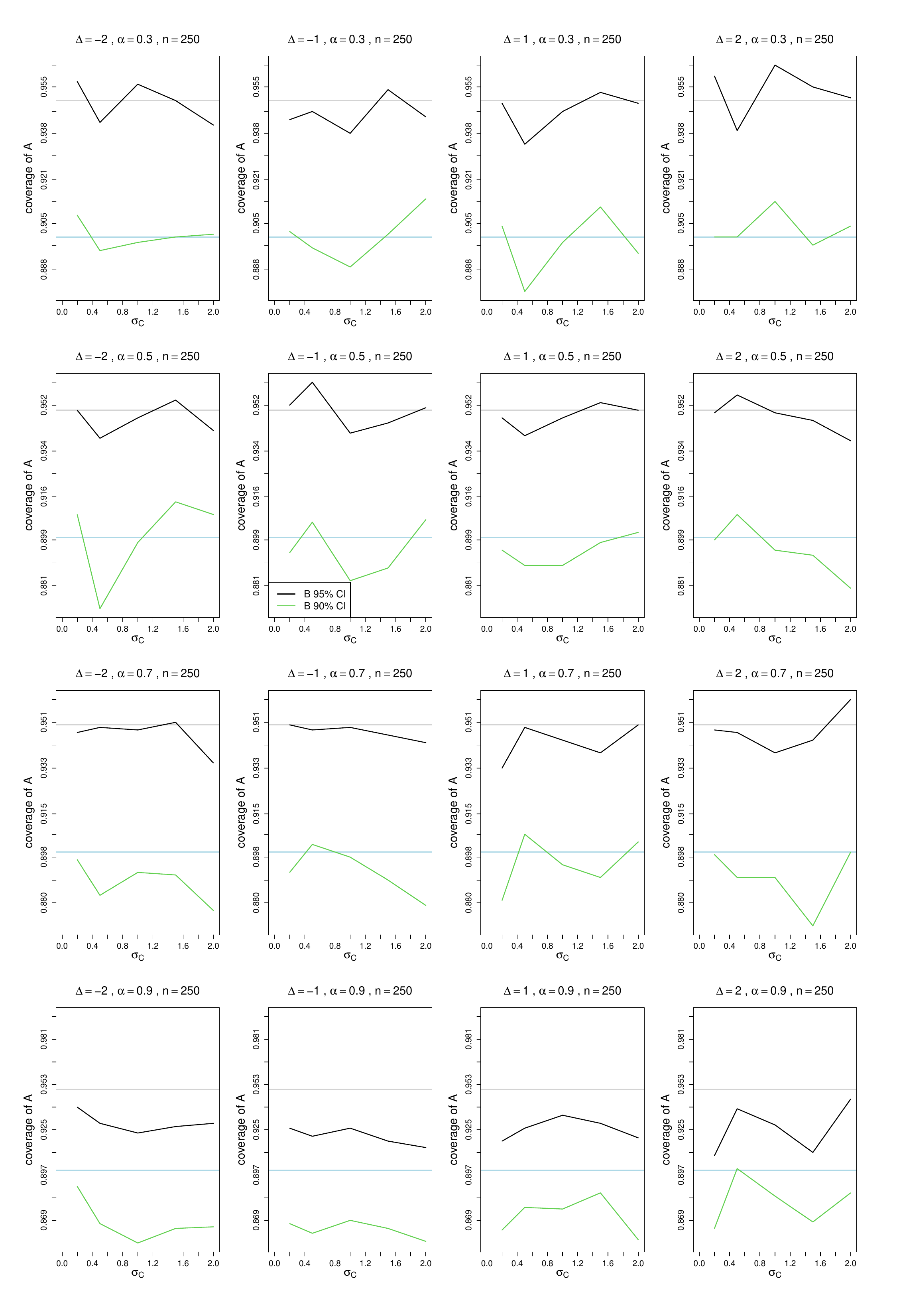}
	\caption{Coverage of $\mathsf{A}$ for $k = k_{max}$ at 95\% (black) and 90\% (green) nominal level as a function of the standard deviation in the Control group $\sigma_C$ when $\mu_C = 1$, $n = 250$, $\alpha = 0.3,\;0.5,\;0.7,\;0.9$, and $\Delta = -2,\;-1,\; 1,\;2$. Light grey line at .95 and light blue line at .90.}
   	\label{FigB4}
\end{figure}

\begin{figure}[ht]
	\centering
	\includegraphics[scale=0.35]{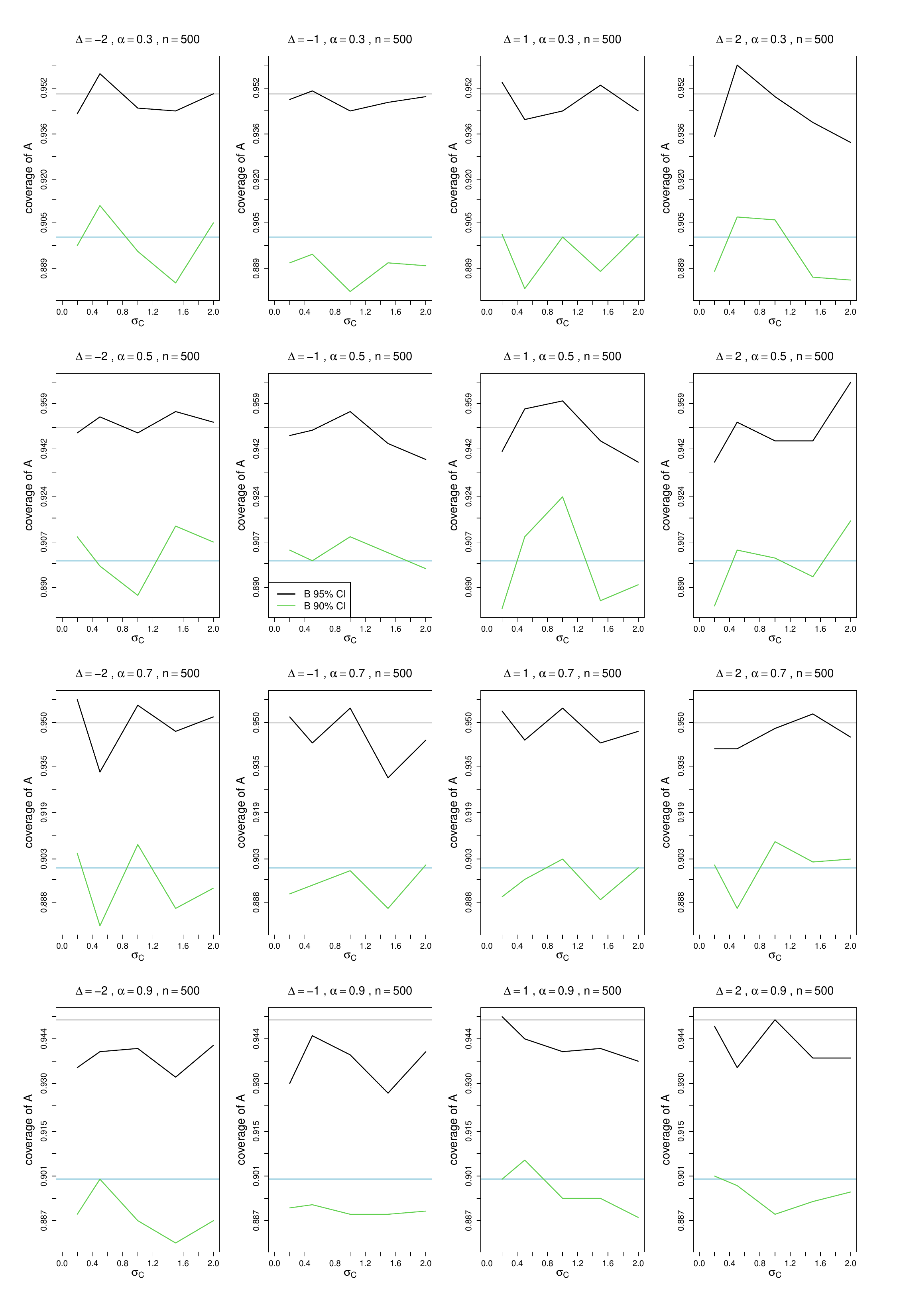}
	\caption{Coverage of $\mathsf{A}$ for $k = k_{max}$ at 95\% (black) and 90\% (green) nominal level as a function of the standard deviation in the Control group $\sigma_C$ when $\mu_C = 1$, $n = 500$, $\alpha = 0.3,\;0.5,\;0.7,\;0.9$ ,and $\Delta = -2,\;-1,\; 1,\;2$. Light grey line at .95 and light blue line at .90.}
   	\label{FigB5}
\end{figure}

\begin{figure}[ht]
	\centering
	\includegraphics[scale=0.35]{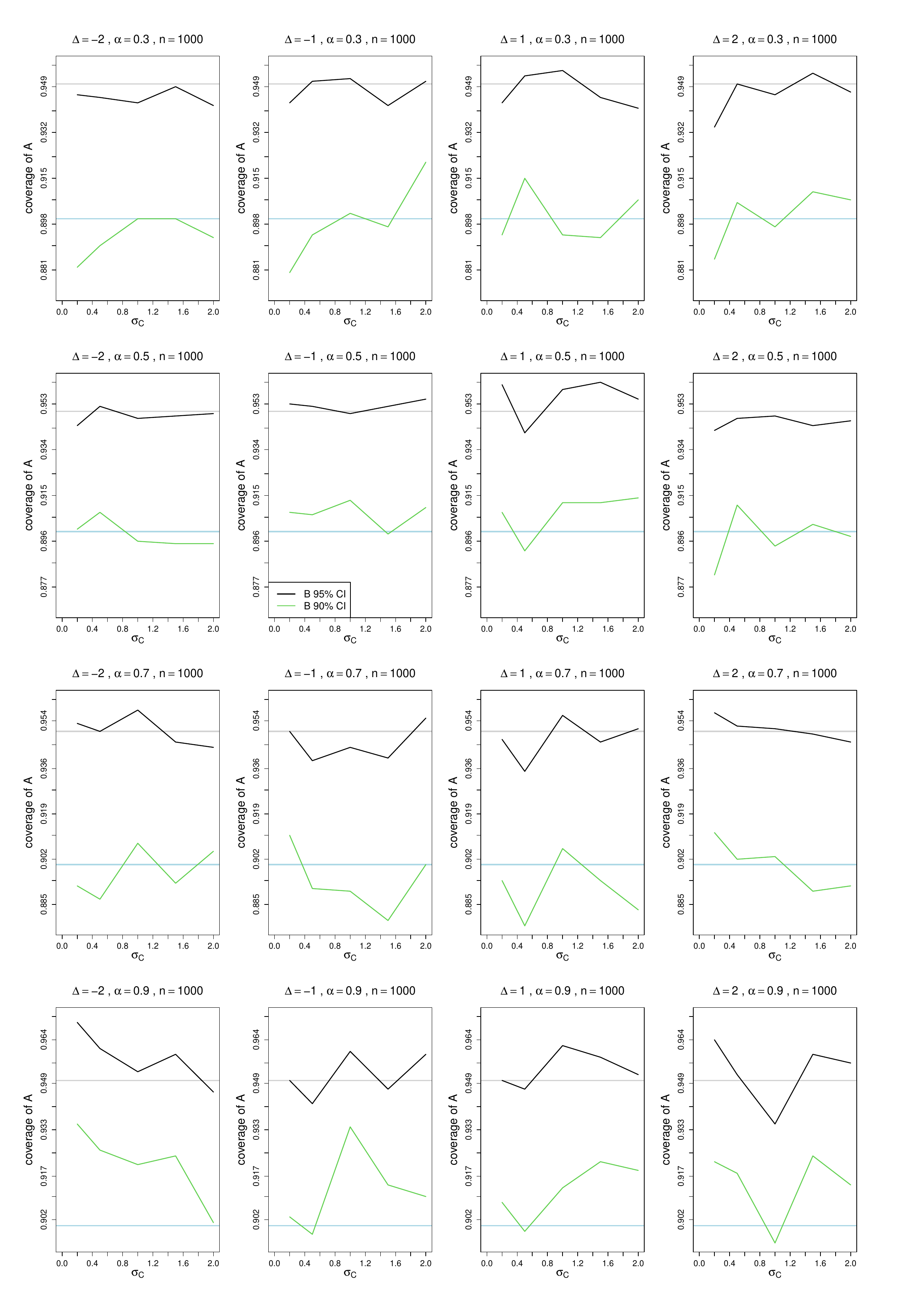}
	\caption{Coverage of $\mathsf{A}$ for $k = k_{max}$ at 95\% (black) and 90\% (green) nominal level as a function of the standard deviation in the Control group $\sigma_C$ when $\mu_C = 1$, $n = 1000$, $\alpha = 0.3,\;0.5,\;0.7,\;0.9$, and $\Delta = -2,\;-1,\; 1,\;2$. Light grey line at .95 and light blue line at .90.}
   	\label{FigB6}
\end{figure}

\clearpage
\setcounter{figure}{0}
\setcounter{section}{0}
\renewcommand{\thefigure}{C.\arabic{figure}}

\section*{Appendix C: Relative width of the bootstrap confidence intervals for $\mathsf{A}$ }

Each figure corresponds to $\mu_C = 1$ and a value of the within-group sample size $n$ (= 20, 50, 100, 250, 500, 1000).\\

\noindent For each combination of a value of $\alpha = \sigma_T / \sigma_C$ (= 0.3, 0.5, 0.7, 0.9) and a value of  the difference in means $\Delta = \mu_T - \mu_C$ (= $-2, -1, 1, 2$), a panel plots, versus the standard deviation in the Control group $\sigma_{C}$  (= $0.2, 0.5, 1, 1.5, 2$), the relative width of the left (black) and right (green) parts of the 90\% (dashed) and 95\% (solid) parametric-bootstrap-based confidence intervals ($B = 2000$) for $\mathsf{A}$.

The relative widths of the left and right parts are based on the lower and upper confidence limits: $100 (\mathsf{A} - \hat{\mathsf{A}}_{\beta / 2}^{\ast}) / | \mathsf{A} |$ and $100 (\hat{\mathsf{A}}_{1 - \beta / 2}^{\ast} - \mathsf{A}) / | \mathsf{A} |$.
\clearpage

\begin{figure}[ht]
	\centering
	\includegraphics[scale=0.35]{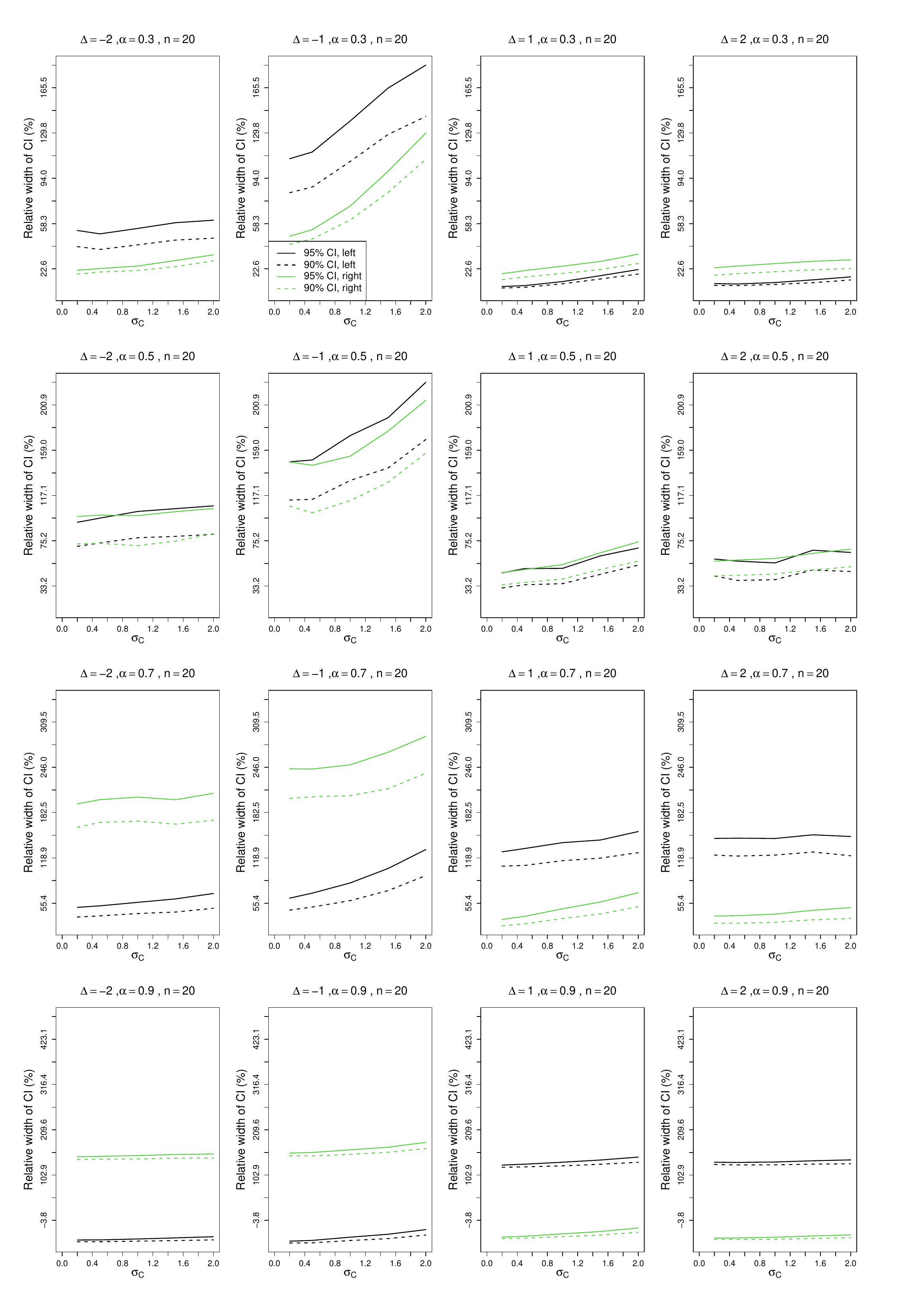}
	\caption{Relative widths of left part (black) and right part (green) of 90\% (dashed) and 95\% (solid) bootstrap confidence intervals for $\mathsf{A}$ for $k = k_{max}$ as a function of the standard deviation in the Control group $\sigma_C$ when $\mu_C = 1$, $n = 20$, $\alpha = 0.3,\;0.5,\;0.7,\;0.9$, and $\Delta = -2,\;-1,\; 1,\;2$. }
   	\label{FigC1}
\end{figure}

\begin{figure}[ht]
	\centering
	\includegraphics[scale=0.35]{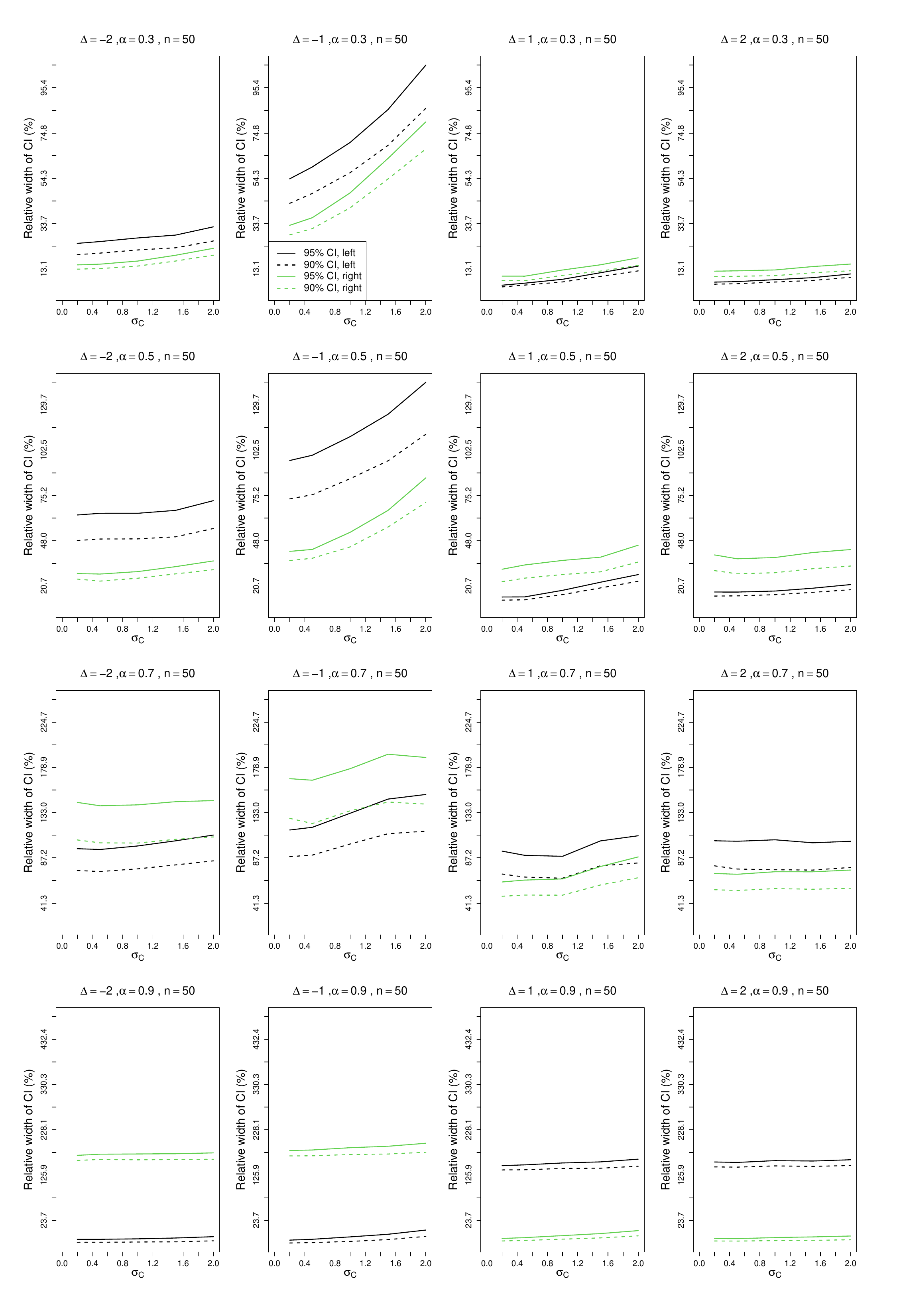}
	\caption{Relative widths of left part (black) and right part (green) of 90\% (dashed) and 95\% (solid) bootstrap confidence intervals for $\mathsf{A}$ for $k = k_{max}$ as a function of the standard deviation in the Control group $\sigma_C$ when $\mu_C = 1$, $n = 50$, $\alpha = 0.3,\;0.5,\;0.7,\;0.9$, and $\Delta = -2,\;-1,\; 1,\;2$. }
   	\label{FigC2}
\end{figure}

\begin{figure}[ht]
	\centering
	\includegraphics[scale=0.35]{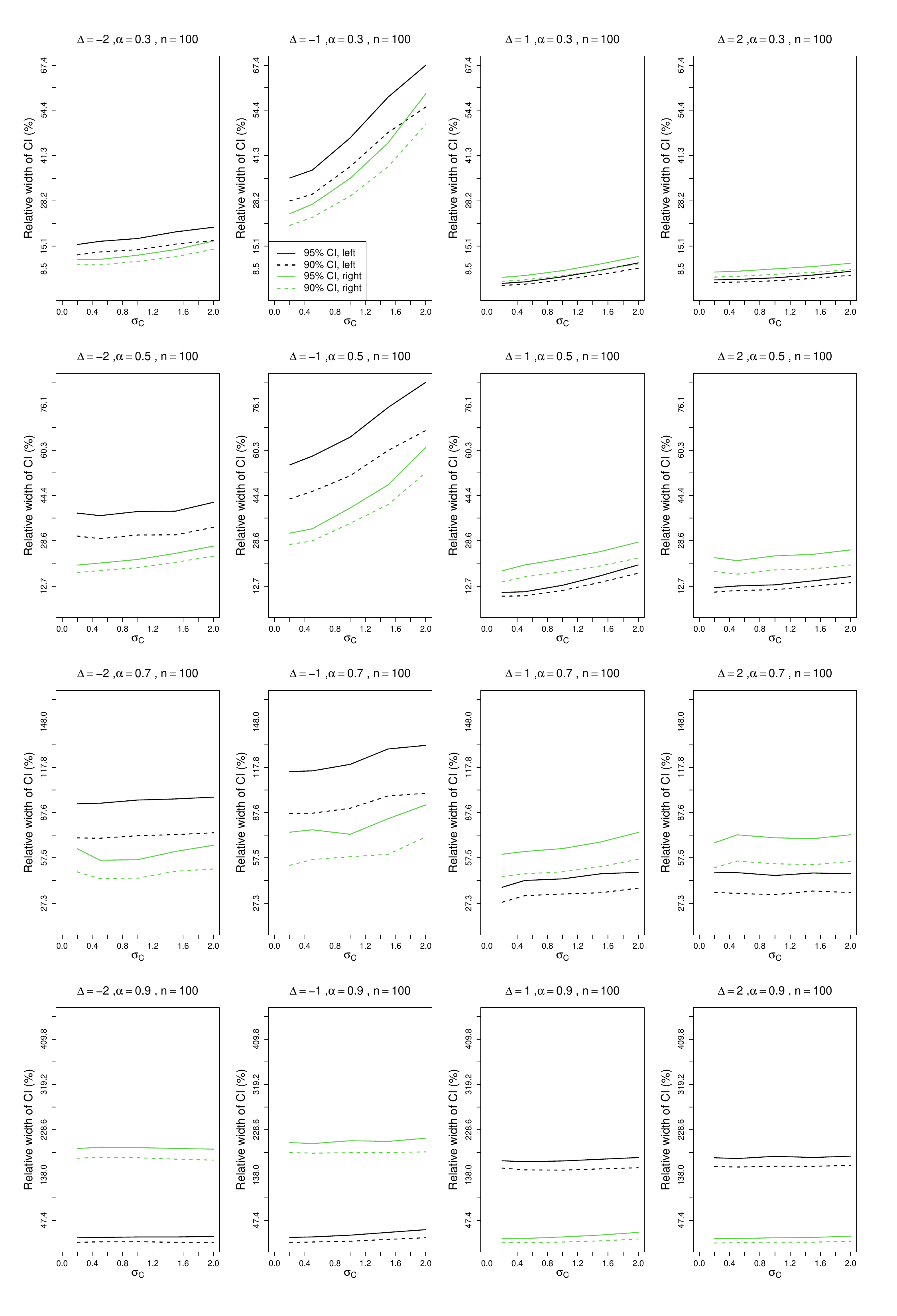}
	\caption{Relative widths of left part (black) and right part (green) of 90\% (dashed) and 95\% (solid) bootstrap confidence intervals for $\mathsf{A}$ for $k = k_{max}$ as a function of the standard deviation in the Control group $\sigma_C$ when $\mu_C = 1$, $n = 100$, $\alpha = 0.3,\;0.5,\;0.7,\;0.9$, and $\Delta = -2,\;-1,\; 1,\;2$. }
   	\label{FigC3}
\end{figure}

\begin{figure}[ht]
	\centering
	\includegraphics[scale=0.35]{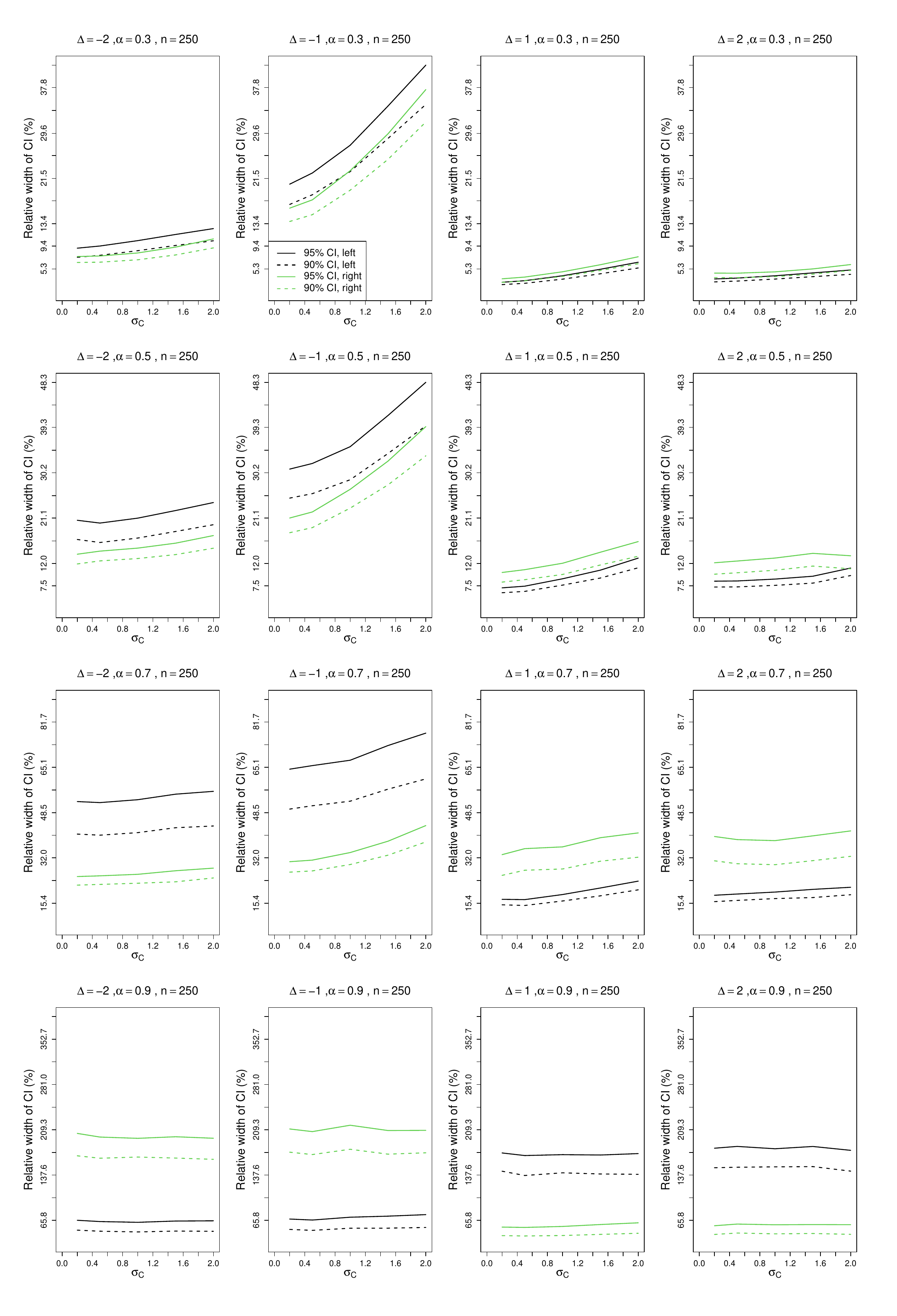}
	\caption{Relative widths of left part (black) and right part (green) of 90\% (dashed) and 95\% (solid) bootstrap confidence intervals for $\mathsf{A}$ for $k = k_{max}$ as a function of the standard deviation in the Control group $\sigma_C$ when $\mu_C = 1$, $n = 250$, $\alpha = 0.3,\;0.5,\;0.7,\;0.9$, and $\Delta = -2,\;-1,\; 1,\;2$. }
   	\label{FigC4}
\end{figure}

\begin{figure}[ht]
	\centering
	\includegraphics[scale=0.35]{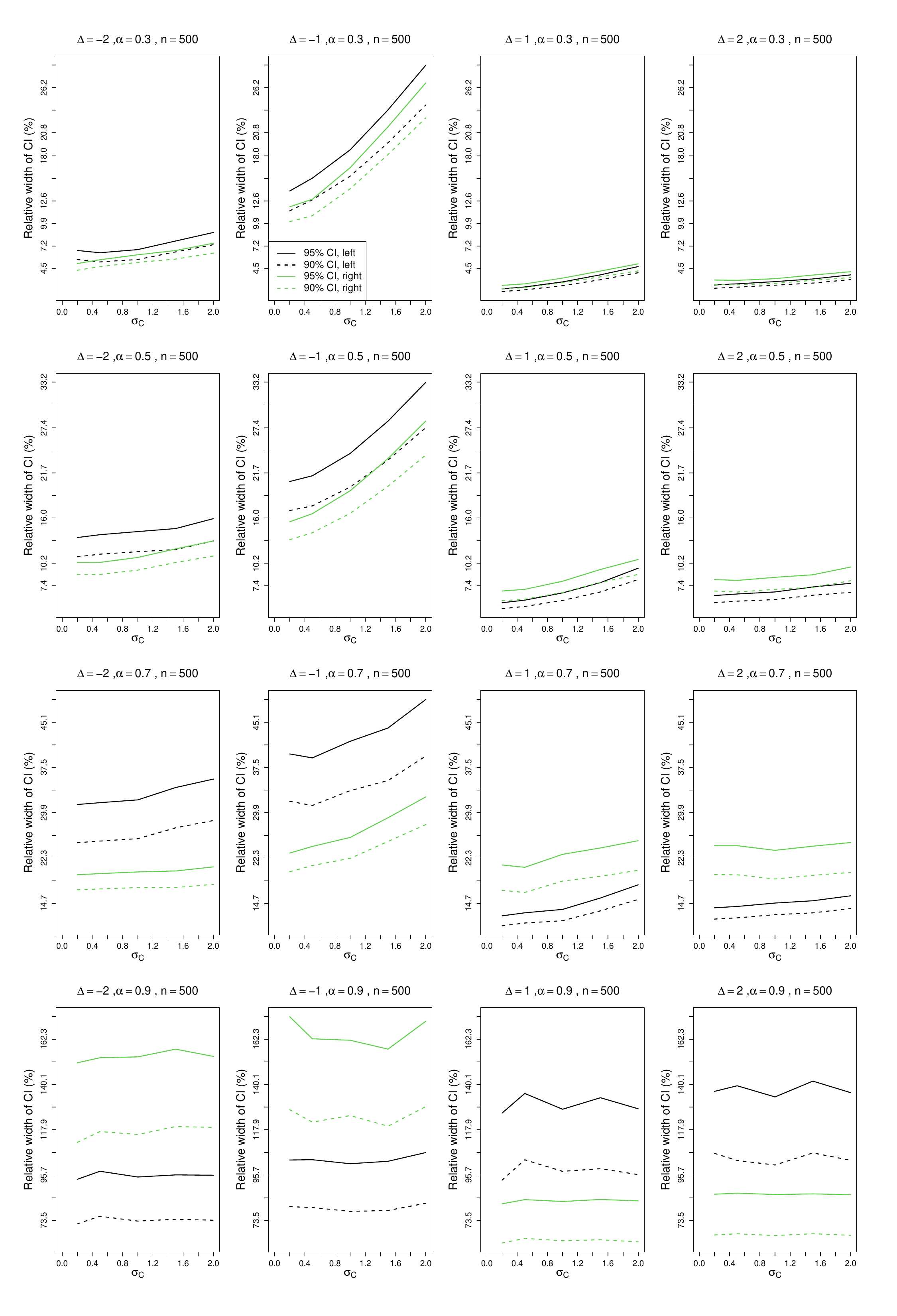}
	\caption{Relative widths of left part (black) and right part (green) of 90\% (dashed) and 95\% (solid) bootstrap confidence intervals for $\mathsf{A}$ for $k = k_{max}$ as a function of the standard deviation in the Control group $\sigma_C$ when $\mu_C = 1$, $n = 500$, $\alpha = 0.3,\;0.5,\;0.7,\;0.9$, and $\Delta = -2,\;-1,\; 1,\;2$. }
   	\label{FigC5}
\end{figure}

\begin{figure}[ht]
	\centering
	\includegraphics[scale=0.35]{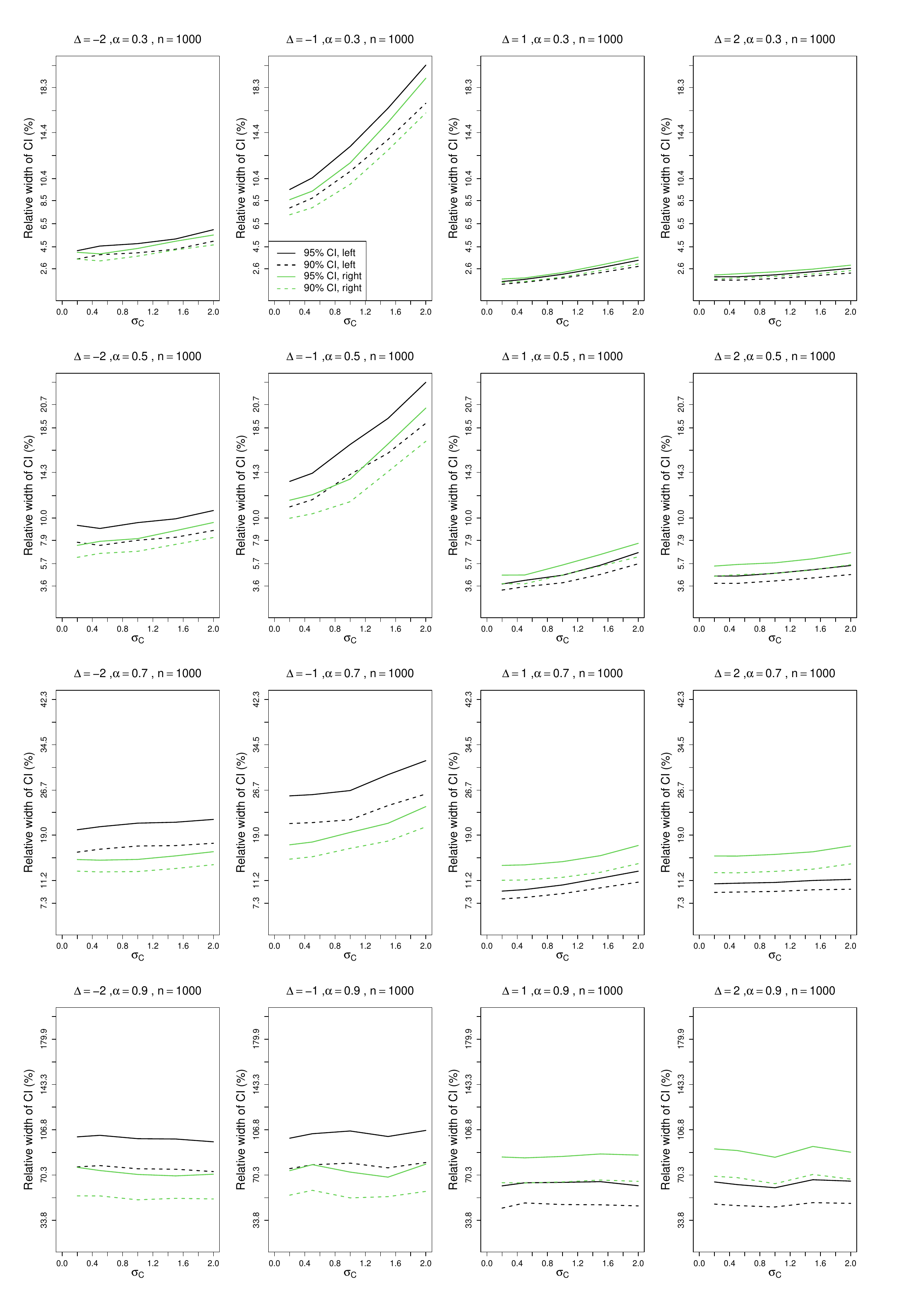}
	\caption{Relative widths of left part (black) and right part (green) of 90\% (dashed) and 95\% (solid) bootstrap confidence intervals for $\mathsf{A}$ for $k = k_{max}$ as a function of the standard deviation in the Control group $\sigma_C$ when $\mu_C = 1$, $n = 1000$, $\alpha = 0.3,\;0.5,\;0.7,\;0.9$, and $\Delta = -2,\;-1,\; 1,\;2$. }
   	\label{FigC6}
\end{figure}
%%%%%%%%%%%%%%%%%%%%%%%%%%%%%%%%%%%%%%%%%%%%%%%%%%%%%%%%%%%%%%%%%%%%%%%%%%%%%%%%%%%%%%%%%%%%%%%%%%%%%%%%%%%%%%%%%%%%%%%%%%%%%%%%%
%%%%%%%%%%%%%%%%%%%%%%%%%%%%%%%%%%%%%%%%%%%%%%% delta_C=2.5 %%%%%%%%%%%%%%%%%%%%%%%%%%%%%%%%%%%%%%%%%%%%%%%%%%%%%%%%%%%%%%%%%%%%%%%%
\clearpage

\end{document}